\documentclass[a4paper,leqno]{article}
\usepackage{amssymb,amsmath,amsfonts,amsthm,
mathrsfs,graphicx}

\newtheorem{theorem}{Theorem}[section]
\newtheorem{corollary}[theorem]{Corollary}
\newtheorem{lemma}[theorem]{Lemma}
\newtheorem{proposition}[theorem]{Proposition}
\newtheorem{definition}[theorem]{Definition}
\newtheorem{asmpt}[theorem]{Assumption}
\theoremstyle{definition}
\newtheorem{remark}[theorem]{Remark}

\newcommand\dslash{d\llap {\raisebox{.75ex}{$\scriptstyle-\!$}}}

\newcommand{\A}{\ensuremath{{\cal A}}}
\newcommand{\B}{\ensuremath{{\cal B}}}

\newcommand{\Con}{\ensuremath{\mathscr{C}}}
\newcommand{\Cinf}{\ensuremath{\mathscr{C}^\infty}}

\newcommand{\D}{\ensuremath{{\mathscr D}}}
\renewcommand{\S}{\mathscr{S}}
\newcommand{\sph}{\mathcal S}

\renewcommand{\P}{\mathfrak{P}}

\newcommand{\mb}[1]{\ensuremath{\mathbb{#1}}}
\newcommand{\N}{\mb{N}}

\newcommand{\R}{\mb{R}}

\newcommand{\G}{\ensuremath{{\cal G}}}
\newcommand{\W}{\ensuremath{{\cal W}}}
\newcommand{\J}{\ensuremath{{\cal J}}}

\newcommand{\WF}{\mathrm{WF}}

\newcommand{\Char}{\ensuremath{\text{Char}}}

\newfont{\bigmath}{cmr12 at 13pt}

\newfont{\grecomath}{cmmi12 at 15pt}

\renewcommand{\d}{\ensuremath{\partial}}

\newfont{\bl}{msbm10 scaled \magstep2}

\newcommand{\beq}{\begin{equation}}
\newcommand{\eeq}{\end{equation}}

\newcommand{\inp}[2]{\langle #1 | #2 \rangle}  
\newcommand{\notmid}{\mid\kern-0.5em\not\kern0.5em}

\renewcommand{\Re}{\ensuremath{\text{Re}}}
\renewcommand{\Im}{\ensuremath{\text{Im}}}

\newcommand{\symb}[1]{\sigma\left\{#1\right\}}

\newcommand{\e}[2]{\langle #1 \rangle^{#2}}

\newcommand{\hf}{\frac{1}{2}}
\newcommand{\wrto}{w.r.t$.$}
\newcommand{\wrt}{w.r.t$.\ $}

\renewcommand{\dslash}{d\hspace{-0.4em}{ }^-\hspace{-0.2em}}
 
\begin{document}

\title{\textsf{\textbf{Fourier-integral-operator approximation of
      solutions to first-order hyperbolic pseudodifferential equations
      I: convergence in Sobolev spaces}}}

\author{J\'{e}r\^{o}me Le Rousseau \\
Laboratoire d'analyse, topologie, probabilit\'{e}s, CNRS UMR 6632 \\
Universit\'{e} Aix-Marseille I, France \\
\texttt{jlerous@cmi.univ-mrs.fr}}
\date{\today}
\maketitle

\begin{abstract} 
  \noindent
  An approximation Ansatz for the operator solution, $U(z',z)$, of a
  hyperbolic first-order pseudodifferential equation, $\d_z +
  a(z,x,D_x)$ with $\Re (a) \geq 0$, is constructed as the composition
  of global Fourier integral operators with complex phases. An
  estimate of the operator norm in $L(H^{(s)},H^{(s)})$ of these
  operators is provided which allows to prove a convergence result for
  the Ansatz to $U(z',z)$ in some Sobolev space as the number of
  operators in the composition goes to $\infty$.\\
  
  \noindent
  \emph{AMS 2000 subject classification: 35L05, 35L80, 35S10, 35S30, 86A15.}  
\end{abstract}
\setcounter{section}{-1}
\newcounter{delta}
\renewcommand{\thedelta}{\arabic{delta}}
\setcounter{delta}{0}

\section{Introduction}
We consider the Cauchy problem
\begin{align}
\label{eq:one-way0}
\d_z u + a(z,x,D_x) u &=0, \ \ \ 0< z\leq Z\\
u \mid_{z=0} &= u_0,
\label{eq:init_cond0}
\end{align}
with $Z>0$ and $a(z,x, \xi)$ continuous with respect to (\wrto) $z$
with values in $S^1(\R^n\times\R^n)$. We denote $D_x = \frac{1}{i}
\d_x$.  Further assumptions will be made on the symbol $a(z,x,\xi)$.
When $a(z,x,\xi)$ is in fact independent of $x$ and $z$ it is natural
to treat such a problem by means of Fourier transformation:
\begin{align*}
  u(z,x') = \iint \exp[i \inp{x'-x}{\xi}
  - z a(\xi)]\ u_0(x) \ \dslash \xi\ d x, 
\end{align*}
where $\dslash \xi := (\frac{1}{2\pi})^n d \xi$.  For this to be well
defined for all $u_0 \in \S(\R^n)$ or some Sobolev space we shall
impose the real part of the principal symbol of $a$ to be
non-negative. When the symbol $a$ depends on both $x$ and $z$ we can
naively expect
\begin{align*}
  u(z,x') \simeq u_1(z,x'):=\iint \exp[i \inp{x'-x}{\xi}
  - z a(0,x',\xi)]\ u_0(x) \ \dslash \xi\ d x
\end{align*}
for $z$ small and hence approximately solve the Cauchy
problem~(\ref{eq:one-way0})-(\ref{eq:init_cond0}) for $z\in
[0,z^{(1)}]$ with $z^{(1)}$ small.  If we want to progress in the 
$z$ direction we have to solve the Cauchy problem
\begin{align*}
\d_z u + a(z,x,D_x) u &=0, \ \ \ z^{(1)}< z\leq Z \\
u(z,.)\mid_{z=z^{(1)}} &= u_1(z^{(1)},.),
\end{align*}  
which we approximatively solve by 
\begin{multline*}
  u(z,x') \simeq u_2(z,x')\\
  :=\iint \exp[i \inp{x'-x}{\xi}
  - (z-z^{(1)}) a(z^{(1)},x',\xi)]\ u_1(z^{(1)},x) \ \dslash \xi\ d x. 
\end{multline*}
Such a procedure then goes on.

If we call $\G_{(z',z)}$ the operator with kernel 
\begin{align*}
  G_{(z',z)}(x',x) = \int \exp[i \inp{x'-x}{\xi} ]
  \exp[ - (z'-z)a(z,x',\xi)]\, \dslash \xi,
\end{align*}
we then see that the procedure described above involves composing such
operators: with chosen values $0\leq z^{(1)}\leq \dots\leq z^{(k)}\leq Z$
we then have 
\begin{align*}
  u_{k+1}(z,x) = \G_{(z,z^{(k)})} \circ \G_{(z^{(k)},z^{(k-1)})} \circ \dots
  \circ\G_{(z^{(1)},0)} (u_0) (x),
\end{align*} 
for $z\geq z^{(k)}$.  We then naturally define for
$\P=\{z^{(0)},z^{(1)},\dots,z^{(N)}\}$, a subdivision of $[0,Z]$ with
$0=z^{(0)}< z^{(1)}<\dots <z^{(N)}=Z$, the operator $\W_{\P,z}$
\begin{align*}
  \W_{\P,z} :=
  \left\{
    \begin{array}{ll}
      \G_{(z,0)} & \text{if } 0\leq z\leq z^{(1)},\\
      \G_{(z,z^{(k)})} {\displaystyle \prod_{i=1}^{k}} 
      \G_{(z^{(i)},z^{(i-1)})} 
      & \text{if } z^{(k)}\leq z\leq z^{(k+1)}.
    \end{array}
  \right.
\end{align*}
The procedure described above yields $\W_{\P,z}(u_0)$ as an
approximation Ansatz for the solution to the Cauchy
problem~(\ref{eq:one-way0})-(\ref{eq:init_cond0}). We denote 
$\Delta_\P=\sup_{i=1,\dots,N} (z_i - z_{i-1})$. The operator
$\G_{(z',z)}$ is often referred to as the thin-slab propagator (see
e.g.~\cite{dHlRW:00,dHlRB:03}).

Note that a similar procedure can be used to show the existence of an
evolution system by approximating it by composition of semigroup
solutions of the Cauchy problem with $z$ 'frozen' in $a(z,x,D_x)$
\cite{Pazy:83,K:70}. Note that the thin-slab propagator $\G_{(z',z)}$
is however not a semigroup nor an evolution family here (see
Section~\ref{sec:3} for simple arguments).

The approximation Ansatz proposed here is a tool to compute
approximations of the exact solution to the Cauchy
problem~(\ref{eq:one-way0})-(\ref{eq:init_cond0}). Such computations
in the case of geophysical problems can be found in~\cite{dHlRW:00}.
In exploration seismology one is confronted with solving equations of
the type
\begin{align}
  \label{eq:geoph1}
  (\d_z -i b(z,x,D_t,D_x) +c(z,x,D_t,D_x)) v &=0, \\ 
  \label{eq:geoph2}
  v(0,.) &=
  v_+(0,.),
\end{align}
where $t$ is time, $z$ is the vertical coordinate and $x$ is the
lateral or transverse coordinate. The operators $b$ and $c$ are of
first order, wiht real principal parts, $b_1$ and $c_1$, and
$c_1(z,x,\tau,\xi)$ is non-negative. We see that the Cauchy
problem~(\ref{eq:one-way0})-(\ref{eq:init_cond0}) studied here is more
general. The Cauchy problem~(\ref{eq:geoph1})-(\ref{eq:geoph2}) is
obtain by a (microlocal) decoupling of the up-going and down-going
wavefields in the acoustic wave equation (see Appendix~\ref{app:A} and
\cite{Stolk:04-bis} for details). The proposed Ansatz can then be a
tool to approximate {\em in practice} the exact solution of the Cauchy
problem~(\ref{eq:geoph1})-(\ref{eq:geoph2}) for the purpose of imaging
the Earth's interior \cite{dHlRW:00,dHlRB:03}. As explained in
Appendix~\ref{app:A} the operator $c$ works as a damping term that
suppresses singularities in the microlocal region where its symbol
does not vanishes. This effect is actually recovered in the proposed
Ansatz.  Geophysists are not only interested in the convergence of
this Ansatz to the exact solution of the Cauchy
problem~(\ref{eq:geoph1})-(\ref{eq:geoph2}) but they expect the
wavefront set of the approximate solution to be close, in some sense,
to that of the exact solution because seismic imaging aims at imaging
the singularities in the subsurface (see for instance
\cite{SdH:02,deHoop:04}).  We shall investigate the microlocal
properties of the proposed Ansatz in Part II, written in collaboration
with G\"unther H\"ormann.

In the present paper, we are interested in the analysis of the
convergence of the approximation scheme $\W_\P$ in Sobolev spaces.
Section~\ref{sec:1} introduces the Cauchy problem we study and the
precise assumptions made on the symbol $a(z,x,\xi)$, especially on the
real part, $c_1$, and imaginary part, $-b_1$, of its principal symbol.
In Section~\ref{sec:2}, we shall at first concentrate our study on the
operator $\G_{(z',z)}$, yet to be properly defined. Under some
assumptions on $a(z,x,\xi)$, we shall prove that $\G_{(z',z)}$ is a
global Fourier integral operator (FIO) with complex phase and that it
maps $\S$ into $\S$, $\S'$ into $\S'$ and $H^{(s)}$ into $H^{(s)}$ for
any $s$. An estimation of $\| \G_{(z',z)} \|_{(H^{(s)},H^{(s)})}$ will
be the first step towards the analysis in Section~\ref{sec:3} of the
convergence of $\W_{\P,z}$. In fact we prove that for $z'-z$
sufficiently small then (Theorem~\ref{theorem:H^s estimate})
\begin{align*}
  \| \G_{(z',z)}\|_{(H^{(s)},H^{(s)})}\leq 1+ |z'-z| M,
\end{align*}
for some constant $M$. Such an estimate is achieved by the
analysis of the behavior of the symbol $\exp[-\Delta c_1]$ as an
  element of $S^0_\hf$, in particular as $\Delta=z'-z$ goes to zero.

In Section~\ref{sec:3} we study the convergence of the Ansatz
 $\W_{\P,z}(u_0)$ to the solution of the Cauchy
 problem~(\ref{eq:one-way0})-(\ref{eq:init_cond0}) in Sobolev
 spaces as $\Delta_\P$ goes to $0$. A convergence in norm of
 $\W_{\P,z}$ to the solution operator of the Cauchy
 problem~(\ref{eq:one-way0})-(\ref{eq:init_cond0}) is actually
 obtained (Theorem~\ref{theorem: convergence result}):
  \begin{align*}
  \lim_{\Delta_\P \to 0} 
  \| \W_{\P,z}  - U(z,0) \|_{(H^{(s+1)},H^{(s)})} = 0, 
  \end{align*}
with a convergence rate of order $\hf$.

 At the end Section~\ref{sec:3} we relax some regularity property of
 the symbol $a(z,.)$ \wrt $z$ by the introduction of another, yet
 natural, Ansatz: following \cite{N.Kumano-go:95}, the thin-slab
 propagator, $\G_{(z',z)}$, is replaced by the operator
 $\widehat{\G}_{(z',z)}$ with kernel
\begin{align*}
  \widehat{G}_{(z',z)}(x',x) = \int \exp[i \inp{x'-x}{\xi} ]
  \exp[ - \textstyle{\int_z^{z'}}a(s,x',\xi) ds]\, \dslash \xi.
\end{align*}

In Part II, we shall focus on the microlocal aspects of the operator
$\W_{\P,z}$ and how it propagates the singularities of the initial
condition $u_0$.  We shall show that the wavefront set of
$\W_{\P,z}(u_0)(z,.)$ converges in some sense to that of the solution
$u(z,.)$ of the Cauchy
problem~(\ref{eq:one-way0})-(\ref{eq:init_cond0}) as $\Delta_\P$ goes
to $0$.

Multi-composition of FIOs to approximate solutions of Cauchy problems
where first proposed in \cite{KuTaTO:78} and \cite{KuTa:79} where the
exact solution operator of a first order hyperbolic system is
approximated with a different Ansatz, up to a regularizing operator.
The technique is based on the computation and the estimation of the
phase functions and amplitudes of the FIO resulting from these
multi-products, a result know as the Kumano-go-Taniguchi theorem.
The technique was then further applied to Schr\"odinger equations
with specific symbols \cite{KiKu:81,N.Kumano-go:95}. In these
latter works the multi-product in also interpreted as an iterated
integral of Feynman's type and convergence is studied in a weak
sense. In \cite{KiKu:81} a convergence result in $L^2$ is proved.
This is the type of results sought here for first order hyperbolic
equations. We however do not use the apparatus of multi-phases and
rather focus on estimating the Sobolev regularity of each term in the
multi-product of FIOs in the proposed Ansatz. While the resulting
product is an FIO, we do not compute its phase and amplitude.

In this paper, when the constant $C$ is used, its value may change
from one line to the other. If we want to keep track of the value of a
constant we shall use another letter. When we shall write that a
function is bounded \wrt $z$ and/or $\Delta$ we shall actually mean
that $z$ is to be taken in the interval $[0,Z]$ and $\Delta$ in some
interval $[0,\Delta_{max}]$ unless otherwise stipulated. We shall
generally write $X$, $X'$, $X''$, $X^{(1)}$, \dots, $X^{(N)}$ for
$\R^{n}$, according to variables, e.g., $x$, $x'$, \dots, $x^{(N)}$.

In the present paper, symbol spaces are spaces of global symbols; a
function $a \in \Cinf(\R^n\times \R^p)$ is in $S^m_{\rho,
  \delta}(\R^n\times \R^p)$, $0<\rho\leq 1$, $0\leq \delta <1$, if for
all multi-indices $\alpha$, $\beta$ there exists $C_{\alpha \beta}>0$
such that
\begin{align*}
  |\d_x^\alpha \d_\xi^\beta  a (x,\xi)| \leq C_{\alpha \beta}\: 
  (1+|\xi|)^{m -\rho |\beta| + \delta |\alpha|}, \ \ x \in \R^n, \ \xi \in \R^p.
\end{align*}
The best possible constants $C_{\alpha \beta}$, i.e., 
\begin{align*}
  p_{\alpha \beta} (a) := \sup_{(x,\xi)\in \R^n\times \R^p}
  (1+|\xi|)^{-m +\rho |\beta| - \delta |\alpha|} |\d_x^\alpha
  \d_\xi^\beta a (x,\xi)|,
\end{align*}
define seminorms that turn $S^m_{\rho, \delta}(R^n\times \R^p)$ into a
Fr\'{e}chet space. As usual we write $S^m_\rho (\R^n\times \R^p)$ in the case
$\rho=1-\delta$, $\hf\leq \rho < 1$,  and $S^m(\R^n\times \R^p)$ in the case
$\rho=1$, $\delta=0$.

We shall use, in a standard way, the notation $\#$ for the composition
of symbols of pseudodifferential operators ($\psi$DO). When given an
amplitude $p(x,y,\xi) \in S^m_{\rho,\delta} (X\times X\times \R^n)$,
$\rho \geq \delta$, we shall also use the notation $\symb{p}(x,\xi)$
for the symbol of the pseudodifferential operator with amplitude
$p$. For $p \in S^m_{\rho,\delta} (X\times \R^n)$ we shall write
$p^\ast$ for the symbol of the adjoint  operator. When composing $\psi$DOs
or computing adjoints of $\psi$DOs we shall make use of the
oscillatory integral representation of the resulting symbol instead of
asymptotic series for two reasons. First, we aim at estimating
operator norm in $L(H^{s},H^{s})$ while using asymptotic series
representations for symbols yields results up to regularizing
operators which operator norms cannot be controlled. Second, we shall
consider symbols in $S^m_\rho$, for some m, including the case
$\rho=\hf$ for which the asymptotic formulae of the calculus of
$\psi$DOs cease to hold.

For $r \in \R$ we let $E^{(r)}$ be the $\psi$DO with
symbol $\e{\xi}{r}:=(1+ |\xi|^2)^{r/2}$. The operator $E^{(r)}$  maps
$H^{(s)}(X)$ onto $H^{(s-r)}(X)$ unitarily for all $s\in \R$ with
$E^{(-r)}$ being the inverse map.


\section{The homogeneous first-order hyperbolic equation}
\label{sec:1}
Let $s \in \R$ and $Z>0$. We consider the Cauchy problem
\begin{align}
\d_z u + a(z,x,D_x) u &=0, \ \ \ 0< z\leq Z,
\label{eq:one-way}\\
u\mid_{z=0} &= u_0 \in H^{(s+1)}(\R^{n}),
\label{eq:init_cond}
\end{align}
where the symbol $a(z,x,\xi)$ satisfies the following assumption
\begin{asmpt}
\label{assumpt:general assumption}
\begin{align*}
a_z (x,\xi) = a(z,x,\xi) = -i\, b(z,x,\xi) + c(z,x,\xi)
\end{align*}
where $b \in \Con^0([0,Z],S^1(\R^{n}\times \R^{n}))$, with real principal
 symbol $b_1$ homogeneous of degree 1 for $|\xi|$ large enough and $c
 \in \Con^0([0,Z],S^1(\R^{n}\times \R^{n}))$ with non-negative principal
 symbol $c_1$ homogeneous of degree 1 for $|\xi|$ large
 enough. Without loss of generality we can assume that $b_1$ and $c_1$
 are homogeneous of degree 1 for $|\xi|\geq 1$.
\end{asmpt}
In Section~\ref{sec:3} we shall further make the following assumption.
\begin{asmpt}
  \label{assumpt:Lipschitz assumption}
  The symbol $a(z,.)$ is assumed to be in ${\mathscr L}([0,Z],S^1(\R^{n}\times
  \R^{n}))$, i.e. Lipschitz continuous \wrt $z$ with values in
  $S^1(\R^{n}\times \R^{n})$, in the sense that,
\begin{align*}
  a(z',x,\xi) - a(z,x,\xi) = (z'-z) \tilde{a}(z',z,x,\xi),\  \ 
  0 \leq z \leq z' \leq Z
\end{align*}
  with $\tilde{a}(z',z,x,\xi)$ bounded \wrt $z'$ and $z$ with values
in $S^1(\R^{n}\times \R^{n})$.
\end{asmpt}
Weaker assumptions will also be formulated in Section~\ref{sec:3}, for
instance by the introduction of another approximating Ansatz.

We denote by $a_1= -i b_1 +c_1$ the principal symbol of $a$ and write
$b=b_1 + b_0$ with $b_0 \in \Con^0([0,Z],S^0(\R^{n}\times \R^{n}))$
and $c=c_1 + c_0$ with $c_0 \in \Con^0([0,Z],S^0(\R^{n}\times
\R^{n}))$.  Assumption \ref{assumpt:general assumption} ensures that
the hypotheses (i)-(iii) of Theorem 23.1.2 of \cite{Hoermander:V3} are
satisfied. Then there exists a unique solution in
$\Con^0([0,Z],H^{(s+1)}(\R^n)) \cap \Con^1([0,Z],H^{(s)}(\R^n))$ to
the Cauchy problem (\ref{eq:one-way})-(\ref{eq:init_cond}).

Furthermore, we have the following energy estimate \cite[Lemma
23.1.1]{Hoermander:V3} for any function in
$\Con^1([0,Z],H^{(s)}(\R^n)) \cap \Con^0([0,Z],H^{(s+1)}(\R^n))$
\begin{multline}
  \label{eq:energy estimate}
  \sup_{z\in[0,Z]} \exp [-\lambda z]\ \| u(z,.)\|_{H^{(s)}}
  \leq \| u(0,.)\|_{H^(s)} \\
  + 2 \int_0^Z \exp[-\lambda z]\ \|\d_z u + a_z(x,D_x)u\|_{H^{(s)}}
  dz,
\end{multline}
with $\lambda$ large enough ($\lambda$ solely depending on $s$).

By Proposition 9.3 in \cite[Chapter VI]{EN:99} the family of operators
$(a_z)_{z\in [0,Z]}$ generates a strongly continuous evolution system.
We denote $U(z',z)$ the corresponding evolution system:
\begin{align*}
  U(z'',z') \circ U(z',z) &= U(z'',z),\ \ Z \geq z''\geq z'\geq z \geq
  0.
 \intertext{Then the Cauchy problem
  (\ref{eq:one-way})-(\ref{eq:init_cond}) reads }\d_z U(z,z_0)u_0 +
  a(z,x,D_x) U(z,z_0)u_0 &=0, \ \ 0\leq z_0<z\leq Z, \\ U(z_0,z_0)u_0
  &= u_0 \in H^{(s+1)}(\R^{n})
\end{align*}
while $U(z,z_0) u_0 \in H^{(s+1)}(\R^{n})$ for all $z\in [z_0,Z]$.



\section{The thin-slab propagator. Regularity properties.}
\label{sec:2}
We follow the terminology introduced in \cite[Sections
25.4-5]{Hoermander:V4} for FIOs with complex phase.  Let $z',z \in
[0,Z]$ with $z'\geq z$ and let $\Delta:= z'-z$. Define $\phi_{(z',z)}
\in \Cinf (X'\times X \times \R^n)$ as
\begin{multline}
  \label{eq: phase function}
  \phi_{(z',z)}(x',x,\xi):= \inp{x'-x}{\xi}
  +i \Delta a_1(z,x',\xi)\\
  = \inp{x'-x}{\xi}  + \Delta b_1(z,x',\xi) +i
  \Delta c_1(z,x',\xi).
\end{multline}
\begin{remark}
  The function $\phi_{(z',z)}$ is assumed to be homogeneous of degree
  1 only when $|\xi| \geq 1$. This however is not an obstacle to the
  subsequent analysis, e.g., FIO properties, since to define such
  operators the phase function need not be homogeneous of degree 1 for
  small $|\xi|$. In the subsequent results concerning the phase
  function and FIOs one will then assume that $|\xi|$ is large enough,
  i.e., $|\xi| \geq 1$.
\end{remark}
\begin{lemma}
  $\phi_{(z',z)}$ is a non-degenerate complex phase function of
  positive type (at any point $(x_0',x_0,\xi_0)$ where $\d_\xi
  \phi_{(z',z)} =0$).
\end{lemma}
\begin{proof} 
  Note that, by Assumption~\ref{assumpt:general assumption},
  $\Im(\phi_{(z',z)}) \geq 0$ and $\phi$ is homogeneous of degree 1;
  $\d_x \phi=0$ implies $\xi=0$. Thus, $\phi$ is a phase function of
  positive type. Inspecting the partial derivatives of $\d_\xi \phi$
  \wrt $x$ we conclude that the differentials
  $d(\d_{\xi_1}\phi), \dots, d(\d_{\xi_{n}}\phi)$ are linearly
  independent.
 \end{proof}

 With $a_0(z,.) \in S^0(X\times\R^{n})$ we have
 $\exp[-\Delta a_0(z,.)] \in S^0(X \times\R^{n})$ by Lemma
 18.1.10 in \cite{Hoermander:V3}. We define
 \begin{align}
   \label{eq: amplitude g}
   g_{(z',z)}(x,\xi) := \exp[-\Delta a_0(z,x,\xi)].
 \end{align}
 We shall keep this notation (for this symbol and others in the
 sequel) but it will be useful however to consider this symbol to
 depend on the parameters $z$ and $\Delta$ instead of $z$ and $z'$ in
 the following analysis. Note that $g_{(z',z)}$ is bounded \wrt $z$
 and $\Cinf$ \wrt $\Delta$ with values in $S^0(X \times\R^{n})$.
 Hence, we may define a distribution kernel $G_{(z',z)}(x',x) \in
 \D'(X'\times X)$
\begin{multline*}
  G_{(z',z)}(x',x) = \int \exp[i \inp{x'-x}{\xi} ]
  \exp[ - \Delta a(z,x',\xi)]\, \dslash \xi\\
  = \int \exp[i \phi_{(z',z)}(x',x,\xi)]\: 
  g_{(z',z)}(x',\xi) \, \dslash \xi
\end{multline*}
as an oscillatory integral. We denote the associated operator by
$\G_{(z',z)}$.  This operator is often referred to as the {\em thin-slab
propagator} (see e.g.~\cite{dHlRW:00,dHlRB:03}). We
show that $\G_{(z',z)}$ is a global FIO in $\R^n$.

Define $\alpha:=(x',x,\xi',\xi)$ and
\begin{align*}
  u_{\theta_j}(\alpha, \theta)
  &= \d_{x_j} \phi_{(z',z)}(x',x,\theta) + \xi_j = - \theta_j + \xi_j,\\
  u_{\xi_j}(\alpha, \theta)
  &= \d_{x_j'} \phi_{(z',z)}(x',x,\theta) - \xi_j'
  = \theta_j - \xi_j' + i \Delta \d_{x_j} a_1(z,x',\theta), \\
  u_{x_j}(\alpha, \theta)
  &= \d_{\theta_j}\phi_{(z',z)}(x',x,\theta)
  = x_j' - x_j +  i \Delta \d_{\xi_j} a_1(z,x',\theta),
\end{align*}
where $j=1,\dots,n$. We denote by $\hat{J}_{(z',z)}$ the ideal in
$\Cinf(\R^{5n})$ generated by the functions $u_{\theta_j},u_{\xi_j}$, and
$u_{x_j}$, and we let $J_{(z',z)}$ be the subset of the functions in
$\hat{J}_{(z',z)}$ that are independent of $\theta$.
\addtocounter{delta}{1}
\begin{lemma}
  \label{lem:generators of  J_(z'z)R}
  There exists $\Delta_\thedelta > 0$, such that, for all $z',z \in
  [0,Z]$, with $z'>z$ and $\Delta=z'-z\leq
  \Delta_\thedelta$, the ideal $J_{(z',z)}$ is generated by the
  functions
  \begin{align}
    v_{\xi_j}(\alpha)
    &= \d_{x_j'} \phi_{(z',z)}(x',x,\xi) - \xi_j'
    \label{eq: generators of J_z'z},\\
    &= \xi_j - \xi_j' + i \Delta \d_{x_j} a_1(z,x',\xi)
    = u_{\xi_j} |_{\theta = \xi}\nonumber,\\
    v_{x_j}(\alpha) &=\d_{\xi_j}\phi_{(z',z)}(x',x,\xi) = x_j' - x_j +
    i \Delta \d_{\xi_j} a_1(z,x',\xi) = u_{x_j} |_{\theta = \xi}
    \nonumber
  \end{align}
  $j=1,\dots,n$.
\end{lemma}
Some of the key arguments of the proof are close to that in the proof
of Theorem 25.4.4 in \cite{Hoermander:V4}.
\begin{proof}
  The ideal $\hat{J}_{(z',z)}$ is also generated by the functions
  \[
  u_{\theta_j}, \ \ \tilde{u}_{\xi_j}:=u_{\theta_j} + u_{\xi_j}= \xi_j - \xi_j' + i
  \Delta \d_{x_j} a_1(z,x',\theta), \ \ u_{x_j},
  \]
  $j=1,\dots,n$. We define $\nu:=(x',\xi',\theta)$, $\mu:=(x,\xi)$. We
  set a point $(\nu_0, \mu_0) = (x_0',\xi_0',\theta_0, x_0,\xi_0)$
  where these generators vanish and we work in a neighborhood of this
  point.  (Note that $\theta_0=\xi_0$.)  Since $z \mapsto a_1(z,.) \in
  S^1(X\times \R^n)$ is bounded we have that $\exists \Delta_\thedelta
  > 0$ such that for $0\leq \Delta \leq \Delta_\thedelta$, and all $z
  \in [0,Z]$,
  \begin{align*}
  \det \d \left(u_{\theta_1},\dots ,u_{\theta_{n}},\tilde{u}_{\xi_1},
    \dots, \tilde{u}_{\xi_{n}},
    u_{x_1},\dots,u_{x_{n}}\right)/\d \nu \neq
  0.
  \end{align*}
  By Theorem 7.5.7 in \cite{Hoermander:V1} we have
  \begin{eqnarray*}
    \left(
      \begin{array}{c}
        x'-x\\
        \xi' -\xi\\
        \theta
      \end{array}
    \right)
    =
    \left(
      \begin{array}{cc}
        Q(\nu,\mu) & P(\nu, \mu)\\
        0 & I_n
      \end{array}
    \right)
    \left(
      \begin{array}{c}
        u_{x}\\
        \tilde{u}_\xi\\
        -u_\theta
      \end{array}
    \right)
    +\left(
      \begin{array}{c}
        \tilde{x}(\mu)\\
        \tilde{\xi}(\mu)\\
        \xi
      \end{array}
    \right)
    \end{eqnarray*}
    where $P$ is a $\Con^\infty$ $2n\times n$ matrix and $Q$ is a
    $\Con^\infty$ $2n\times 2n$ matrix and the functions $\tilde{x}$,
    $\tilde{\xi}$ are also $\Con^\infty$ in a neighborhood of
    $(\nu_0,\mu_0)$. As the functions $w_x(\nu,\mu)
    :=x'-x-\tilde{x}(\mu)$, $w_\xi(\nu,\mu)
    :=\xi'-\xi-\tilde{\xi}(\mu)$, $w_\theta(\nu,\mu) := \theta - \xi$
    have linearly independent differentials, Lemma 7.5.8 in
    \cite{Hoermander:V1} proves that they generate $\hat{J}_{(z',z)}$
    and the proof of that lemma shows that $Q$ is invertible in a
    neighborhood of $(\nu_0,\mu_0)$.  Letting $\theta=\xi$ we have
    \begin{eqnarray*}
      \hspace*{1cm}
      Q(x',x,\xi,\mu)^{-1}\ \left(
      \begin{array}{c}
        w_x(\nu,\mu)\\
        w_\xi(\nu,\mu)
        \end{array}
    \right)
    =
    \left(
      \begin{array}{c}
        u_{x}(x',x,\xi)\\
       \tilde{u}_{\xi}(x',x,\xi)
      \end{array}
    \right)
    =\left(
      \begin{array}{c}
        v_x(\alpha)\\
        v_\xi(\alpha)
      \end{array}
    \right).
    \end{eqnarray*}
    We thus obtained that $\hat{J}_{(z',z)}$ is generated by the functions
    $u_{\theta_j}$, $v_{x_j}$, $v_{\xi_j}$,
    $j=1,\dots,n$.  We then see that $J_{(z',z)}$ is generated by
    $v_{x_j}$, $v_{\xi_j}$, $j=1,\dots,n$.
    In fact, using Theorem 7.5.7 in
    \cite{Hoermander:V1} again, any $\Con^\infty$ function $h(\alpha)$ can
    be locally written in the form
    \begin{align*}
    h(\alpha) = \sum_{1\leq i \leq n} (h_{x_j}(\alpha',\mu) v_{x_j}(\alpha',\mu)
     + h_{\xi_j}(\alpha',\mu) v_{\xi_j}(\alpha',\mu))+ r(\mu),
    \end{align*}
    with $\alpha'=(x',\xi')$ provided that $0\leq \Delta\leq \Delta_\thedelta$.  If $h
    \in J_{(z',z)}$ then $r\in J_{(z',z)}$ and Lemma 7.5.10 in
    \cite{Hoermander:V1} implies that $\forall N \in \N$, $\exists C_N>0$:
    \begin{align*}
    r(\mu) \leq C_N \max (|\Im\ \tilde{x}(\mu)|,
    |\Im\ \tilde{\xi}(\mu)|)^N,
    \end{align*}
    locally. Therefore, Theorem 7.5.12 in \cite{Hoermander:V1} yields $r \in
    I(w_{x},w_{\xi}) = I(v_{x},
    v_{\xi})$; which in turn implies $g \in I(v_{x},
    v_{\xi})$ and thereby completes the proof.
\end{proof}

As the Poisson brackets (for the symplectic $2$-form $\sigma' -
\sigma$ on $T^\ast (X'\times X)$, where $\sigma'$ and $\sigma$ are the
symplectic $2$-forms on $T^\ast(X')$ and $T^\ast(X)$ respectively) of
any two of the functions in (\ref{eq: generators of J_z'z}) vanish
identically we obtain that the ideal generated by these functions is
globally a conic canonical ideal in the sense of \cite[Definition
25.4.1. and Section 25.5]{Hoermander:V4}. The phase function
$\phi_{(z',z)}$ thus defines $J_{(z',z)}$ in the neighborhood of any
point of $J_{(z',z)\R}$: it thus {\em globally} defines $J_{(z',z)}$,
which is then of positive type.  Therefore the operator $\G_{(z',z)}$
is a global FIO with complex phase (see Definitions 25.4.9. and
25.5.1. in \cite{Hoermander:V4}).
\begin{proposition}
  There exists $\Delta_\thedelta>0$ such that if $0\leq
  \Delta=z'-z\leq \Delta_\thedelta$ then the operator $\G_{(z',z)}$ is a
  global Fourier integral operator with complex phase and $G_{(z',z)}
  \in I^0(X'\times X, (J_{(z',z)})', \Omega^{1/2}_{X'\times X})$.
\end{proposition}
We denote the half density bundle on $X'\times X$ by
$\Omega^{1/2}_{X'\times X}$.  Note that $(J_{(z',z)})'$ stands for the
twisted canonical ideal, i.e.  a Lagrangian ideal (see Section 25.5 in
\cite{Hoermander:V4}).

Note that, with the following analysis, we could also consider
$\G_{(z',z)}$ as a global FIO with real phase with amplitude in
$S^0_\hf (X'\times X \times \R^n)$ (see e.g. \cite{Ruzhansky:00}). However
such consideration would be rather technical as one usually restricts
oneself to the type $S^m_\rho$ with $\rho > \hf$ for FIOs (see the
remark at the end of Section 25.1 in \cite{Hoermander:V4}; see Also
\cite[pages 391-392]{MeSj:76}). Viewing the thin-slab propagator
$\G_{(z',z)}$ as a FIO with complex phase is also a good framework to
understand the propagation of singularities in Part II. We shall
however make this interpretation for $\G_{(z',z)}$ in
Proposition~\ref{prop: Hs cont any symbol}, below, to apply a result
of Kumano-go~\cite[Theorem 2.5]{Kumano-go:76}.

We now establish some global continuity properties of the operator
$\G_{(z',z)}$ stated in a slightly more general form (for similar
results with global symbols see for instance \cite{Kumano-go:76},
where phase functions are real and other conditions are imposed on the
phase function).
\begin{lemma}
\label{lemma:S cont}
Let $A$ be an FIO with a kernel of the form
\begin{align*}
K_A(x,y) = \int \exp[i \varphi (x, \xi) - i \inp{y}{\xi}]\, \sigma_A(x,\xi)
\dslash \xi \in \D'(\R^n\times\R^n),
\end{align*}
where $\sigma_A \in S^m(\R^n\times \R^n)$ and $\varphi \in \Cinf(\R^n\times
\R^n)$ is such that $\Im (\varphi(x,\xi)) \geq 0$ and $\varphi$ is
homogeneous of degree 1 in $\xi$, for $|\xi|$ large enough, and $\d_{x_i}
\varphi \in S^1(\R^n\times \R^n)$. Furthermore, for all $i=1,\dots,n$ we assume
$\d_{\xi_i}\varphi(x,\xi) = x_i + f_i(x,\xi)$ where $f_i \in S^0(\R^n\times
\R^n)$. Then A maps $\S$ into $\S$ continuously.
\end{lemma}
\begin{proof}
Let $u \in \S$. We then have
\begin{align*}
|A u(x)| &\leq \int | \sigma_A(x,\xi) (1+ |\xi|)^{-m}| | (1+ |\xi|)^{m} \hat{u}(\xi)|
\,\dslash(\xi)\\
& \leq C \sup_{\xi \in \R^n} | \sigma_A(x,\xi) (1+ |\xi|)^{-m}| \sup_{\xi \in
\R^n} | (1+ |\xi|)^{m+n+1} \hat{u}(\xi)|,
\end{align*}
where $ C = \int (1+|\xi|)^{-n-1} \dslash \xi$. The operator $A$ is hence well
defined from $\S$ into $\Con(\R^n)$. If we differentiate we obtain
\begin{align*}
D_{x_i} A u(x) = \int \exp[i \varphi(x,\xi)] \left( \d_{x_i} \varphi (x,\xi)
\sigma_A(x,\xi) - i \d_{x_i} \sigma_A(x,\xi) \right) \hat{u}(\xi)\, \dslash
\xi\, .
\end{align*}
Noting that $\d_{x_i} \varphi (x,\xi) \sigma_A(x,\xi) - i \d_{x_i} \sigma_A(x,\xi)
\in S^{m+1}(\R^n\times \R^n)$ we similarly have
\begin{align*}
|D_{x_i}A u(x)| &\leq C \sup_{\xi \in \R^n} | (1+ |\xi|)^{m+n+2} \hat{u}(\xi)|\\
& \leq C' \sup_{x \in \R^n} |x^\alpha D^\beta_x u(x)| \, \text{ for some
}\alpha, \beta \geq 0.
\end{align*}
Iterating we find that $Au \in \Cinf(\R^n)$. Integrating by parts we also have
\begin{multline*}
 A(x_j u)(x) = \int \exp[i \varphi(x,\xi)] \left(
\d_{\xi_i} \varphi (x,\xi) \sigma_A(x,\xi) - i \d_{\xi_i} \sigma_A(x,\xi)
\right)
\hat{u}(\xi)\, \dslash(\xi)\\
= x_j A u(x) + \int \exp[i \varphi(x,\xi)] \left( f_i(x,\xi) \sigma_A(x,\xi) - i \d_{\xi_i} \sigma_A(x,\xi)
\right) \hat{u}(\xi)\, \dslash(\xi).
\end{multline*}
Since $f_i(x,\xi) \sigma_A(x,\xi) - i \d_{\xi_i} \sigma_A(x,\xi) \in
S^m(\R^n\times \R^n)$ we obtain
\begin{align*}
  |x_j A u(x)| \leq C \sup_{x \in \R^n} |x^\alpha D^\beta_x u(x)| +
  C\sup_{x \in \R^n} |x^{\alpha'} D^{\beta'}_x u(x)|,
\end{align*}
for some $\alpha, \alpha', \beta, \beta' \geq 0$. Similar estimates
hold for $|x^\alpha D_x^\beta A u(x)|$ because of the hypothesis made
on $f_i$, $i=1,\dots,n$. The operator $A$ thus maps $\S$ into $\S$
continuously.
\end{proof}
To show continuity from $\S'$ into $\S'$ we shall need the following
lemma.
\begin{lemma}
  \label{lemma:7.7.1}
  Let $j,k$ non-negative integers, $u \in \S(\R^n)$, $f \in
  C^{k+1}(\R^n)$ such that
  \begin{align*}
    0\leq \Im\: f (x)\leq C_0,\ x\in \R^n, \  \ \  
    |f^{(r)}(x)| \leq C_r, \ x\in \R^n,  \ 1\leq r\leq k+1.
  \end{align*}
Then we have 
\begin{multline}
  \label{eq:7.7.1}
  \omega^{j+k} \left| \int u(x) (\Im\: f(x))^j \exp[i \omega f(x)]\ d x \right| \\
  \leq C \sum_{|\alpha| \leq k} \sup_{x \in \R^n} |D^\alpha u(x)|
  (|f'(x)|^2 + \Im\: f(x))^{|\alpha|/2 -k},\ \ \omega>0,
\end{multline}
where the constant $C$ is bounded when the function $f$ stays in a
domain of $\Con^{k+1}(\R^n)$ where $C_0$, $C_1, \dots, C_{k+1}$ can be
chosen bounded.
\end{lemma}
\begin{proof}
  The proof is the same as that of Theorem 7.7.1 in
  \cite{Hoermander:V1} where $u \in \Con^k_0(\R^n)$. In fact the
  further assumptions on $f$ made here allow to give global bounds
  that are needed since $u \in \S$ in the present case.
\end{proof}
\begin{lemma}
Let $A$ be an FIO with a kernel of the form:
\begin{align*}
K_A(x,y) = \int \exp[i \inp{x-y}{\xi}+i \gamma (x,\xi)]\, \sigma_A(x,\xi) \dslash(\xi)
\in \D'(\R^n\times\R^n),
\end{align*}
where $\sigma_A \in S^m(\R^n\times \R^n)$ and $\gamma \in S^1(\R^n\times
\R^n)$ is such that $\Im (\gamma(x,\xi)) \geq 0$, and $\gamma$ is homogeneous
of degree 1 in $\xi$, for $|\xi|$ large enough.  Furthermore, we assume that
 there exists $d\geq 0$ such that
\begin{eqnarray}
|\Re\left(\d_x \gamma(x,\xi)\right)| \leq d<1, \ x \in \R^n, \ 
\xi \in \R^n, \ |\xi|=1.
\label{eq:hyp phase}
\end{eqnarray}
Then $A$ maps $\S'$ into $\S'$ continuously. \label{lemma:S' cont}
\end{lemma}
Observe that the differential of $\phi(x,\xi):= \inp{x-y}{\xi}+ \gamma
(x,\xi)$ does not vanish in $\R^{2n}\times \R^n\backslash 0$. The
function $\phi$ is thus a complex phase function. The differentials
$d(\d_{\xi_1}\phi),\dots,d(\d_{\xi_n}\phi)$ are linearly independent.
Hence $\phi$ is a non degenerate complex phase function of positive
type.  Note that by~(\ref{eq:hyp phase}) the function $\inp{x-y}{\xi}
+ \gamma(x,\xi)$ is an operator phase function in the sense of
\cite[Definition 1.4.4.]{Hoermander:71}.
\begin{proof}
  Without loss of generality we may assume that $\gamma$ is
  homogeneous of degree 1 for $|\xi|\geq 1$.  Let $A^t$ be the transpose
  of $A$ and let $u \in \S$, then.
  \begin{align*}
    A^t u (x) & = \int \exp[-i\inp{x}{\xi}]
    \int\exp[i \inp{y}{\xi} + i \gamma(y,\xi)]\,
    \sigma_A (y,\xi)\, u(y)\, d y \dslash \xi    \nonumber
  \end{align*}
  Define 
  \begin{align*}
    v(\xi,\eta) = \int \exp[i \inp{y}{\xi} + i \gamma(y,\xi)]\,
    \sigma_A (y,\eta)\, u(y)\, d y\, ,
    \nonumber
  \end{align*}
  and put $w(\xi) = v(\xi,\xi)$. As $u \in \S$ then $v$ and $w$ are
  both $\Cinf$. Then $A^t u$ is the Fourier transform of $w$. The
  lemma is proved if we show that $u \mapsto w(\xi)$ is continuous
  from $\S$ to $\S$.
  
  Let $\omega = |\xi|\geq 1$ and $\xi_0 = \xi/|\xi| \in \sph^{n-1}$.
  We then have $\inp{y}{\xi} + \gamma(y,\xi) = \omega f(y,\xi)$ with
  $f$ homogeneous of degree 0 in $\xi$, for $|\xi|\geq 1$. Note that $\d_y
  f(y,\xi) = \xi_0 + \d_y \gamma (y,\xi_0)$.  With the assumption made
  on $\d_y \gamma$ we have $|\d_y f(y,\xi)| \geq c >0$.  Applying
  Lemma~\ref{lemma:7.7.1} and estimate (\ref{eq:7.7.1}) we obtain
  \begin{align*}
    \omega^k |v(\xi,\eta)|
    &\leq K_k \sum_{|\alpha|\leq k}
    \sup | D^\alpha_y (\sigma_A(y,\eta) u(y))| \\
    & \leq K_k'(1 + |\eta|)^m
    \sup_{{|\alpha| \leq k} \atop {y \in \R^n}} | D^\alpha u (y)|, \ 
    \omega\geq 1
  \end{align*}
  where the constants $K_k$, $K_k'$ can be chosen uniformly \wrt
  $\xi$, $|\xi| \geq 1$ since the constants $C_0$, $C_1, \dots,
  C_{k+1}$ of Lemma~\ref{lemma:7.7.1} can be chosen bounded (as $\xi_0
  \in \sph^{n-1}$). Now setting $\eta=\xi$ we obtain that for all $k
  \in \N$, $\exists K_k''>0$
  \begin{align}
    (1+|\xi|)^{k-m} | w(\xi) |
    \leq K_k'' \!\sup_{{|\alpha| \leq k} \atop {y \in \R^n}}
    |D^\alpha u (y)|\, , \ \xi\in \R^n, \ |\xi|\geq 1.
    \label{eq: S estimate}
  \end{align}
  We now consider
  \begin{multline*}
    D_{\xi_i} w(\xi) = \int \exp[i \inp{y}{\xi} + i \gamma(y,\xi)]\\
    \big( (y_i + \d_{\xi_i} \gamma(y,\xi)) \sigma_A(y,\xi) - i \d_{\xi_i} \sigma_A(y,\xi)\big)
    u(y) d y\, .
  \end{multline*}
  As $y_i u(y) \in \S$ and $ \d_{\xi_i} \gamma(y,\xi)$ is homogeneous of
  degree 0 for $|\xi|\geq 1$ estimates similar to those in (\ref{eq: S
    estimate}) are valid.
\end{proof}

It is immediate from the structure of $\phi_{(z',z)}$ in (\ref{eq:
  phase function}) that Lemma \ref{lemma:S cont} applies to
$\G_{(z',z)}$. If $\Delta=z'-z$ is small enough we have $|\Delta
\partial_{x_i} b_1(z,x',\xi)|\leq d < 1$, due to
Assumption~\ref{assumpt:general assumption}, and then Lemma
\ref{lemma:S' cont} applies. We thus have
\addtocounter{delta}{1}
\begin{proposition}
  There exists $\Delta_\thedelta>0$ such that if $z',z \in [0,Z]$ with
  $0\leq \Delta:= z'-z \leq \Delta_\thedelta$ then $\G_{(z',z)}$ maps
  $\S$ into $\S$ and $\S'$ into $\S'$ continuously.
\end{proposition}
\begin{remark}
  By the above result,  composition of the two FIOs $\G_{(z'',z')}$ and
  $\G_{(z',z)}$ is thus natural without further requirement such as
  having the operators properly supported.
\end{remark}

We now turn to {\em global} $L^2$ and Sobolev space continuity for the
operator $\G_{(z',z)}$.  We shall use the following lemma.
\addtocounter{delta}{1}
\begin{lemma}
  \label{lemma: change of variable - symbol}
  Let $p_s(x,\xi)$ be bounded \wrt the parameter $s$
  with values in $S^m_\rho(X\times \R^n)$ and
  define
  \begin{align*}
  \tilde{\xi}(\Delta,x,\xi):= \xi - \Delta f(x,\xi)
  \end{align*}
  where $f$ is in $S^1 (X\times\R^n,\R^n)$ and homogeneous of
  degree $1$ in $\xi$, for  $|\xi|\geq 1$.  Then
  \begin{align*}
  \tilde{p}_s(\Delta,x,\xi):=
  p_s (x,\tilde{\xi}(\Delta,x,\xi))
  \end{align*}
  is bounded \wrt $s$ with values in
  $\Cinf([0,\Delta_\thedelta], S^m_\rho(X\times
  \R^n))$ for $\Delta_\thedelta$ small enough.
\end{lemma}
\begin{proof}
  Let $\Delta_\thedelta$ be small enough such that $|\xi - \Delta
  f(x,\xi)| \geq C_0>0$ if $|\xi|=1$ and
  $\Delta\in [0,\Delta_\thedelta]$. We then have
  \begin{align*}
    1+ C_0 |\xi| \leq 1+ |\xi - \Delta f(x,\xi)|\leq 1+ C_1|\xi|, \ 
    \xi \in \R^n, \ |\xi|\geq 1, \ \Delta\in [0,\Delta_\thedelta].
  \end{align*}
  This inequality yields the proper estimates for $\d_{x}^\alpha
  \d_{\xi}^\beta \tilde{p}_s$ to prove that $\tilde{p}_s \in
  S^m_\rho(X\times\R^n)$. As derivatives \wrt $\Delta$ do not
  affect the symbol order and type, the proof is finished. Bounds
  \wrt the parameter $s$ come naturally. The proof is complete.
\end{proof}

Following \cite{Stolk:04} we introduce
\begin{definition}
  \label{definition:special damping}
  Let $L\geq 2$. A symbol $q(z,.)$ bounded \wrt $z$ with
  values in $S^1(\R^{p}\times \R^{r})$ is said to
  satisfy Property (\ref{eq:P_L}) if it is non-negative and satisfies
  \begin{multline}
    \label{eq:P_L}
    \tag{$P_L$}
    | \d^\alpha_{y} \d_{\eta}^\beta q(z,y,\eta)|
    \leq
    C (1+|\eta|)^{-|\beta| + (|\alpha| + |\beta|)/L}\\
    (1+ q(z,y,\eta))^{1 - (|\alpha| + |\beta|)/L},\ \
    z \in [0,Z],\ y \in \R^{p},\ \eta \in \R^r.
  \end{multline}
  We then set $\rho=1-1/L$ and $\delta=1/L$. 
\end{definition}
\begin{remark}
\label{rem: P_L for alpha + beta >= L}
Suppose $q(z,.)$ as in Definition~\ref{definition:special damping} and
$|\alpha| + |\beta| \geq L$ then 
\begin{align*}
  (1+ |\eta|)^{1 - (|\alpha| + |\beta|)/L} \leq C(1+ q(z,y,\eta))^{1 -
  (|\alpha| + |\beta|)/L}, z \in [0,Z],\ y \in \R^{p},\ \eta \in \R^r.
\end{align*}
Estimate~(\ref{eq:P_L}) is thus clear in this case.
\end{remark}
Examples of symbols with such a property with $L>2$ are given in
\cite{Stolk:04}. In fact we prove that $c_1$ satisfies Property
(\ref{eq:P_L}) for $L=2$.
\begin{lemma}
  \label{lemma: PL with L=2 for c_1}
  Let $q(z,y,\eta)$ be bounded \wrt $z$ with values in
  $S^1(\R^{p}\times \R^{r})$.  If $q\geq 0$ then $q$ satisfies
  Property (\ref{eq:P_L}) for $L=2$.
\end{lemma}
\begin{proof}
  Bounds \wrt $z$ are natural; we shall omit the dependence on $z$ in
  the proof for concision.  We have to prove that
  \begin{align*}
    |\d_y^\alpha \d_{\eta}^\beta q|
    \leq C\: (1+ |\eta|)^{\hf(|\alpha|-|\beta|)}\
    (1+q)^{1 -\hf(|\alpha|+|\beta|) }
  \end{align*}
  The property is clearly true for $|\alpha|+|\beta|=0$ and for
  $|\alpha|+|\beta|\geq 2$ by the remark above. Let us now treat the
  case $|\alpha|+|\beta|= 1$.  For this we use Landau's inequality:
  let $f\in C^2(\R)$ with $f\geq 0$ and $f''$ is bounded then (see
  \cite[page 40]{Duistermaat:96} and \cite[Lemma
  7.7.2]{Hoermander:V1})
  \begin{align*}
    |f'(t)|\leq 2 \: (f(t))^{\hf}\ \left(\sup_{t \in \R} |f''(t)|\right)^{\hf}.
  \end{align*}
  We first treat the case $|\alpha|=1$.
  Define $p(y,\eta) = (1+ |\eta|^2)^{-\hf}\ q(y,\eta)$. Then $p
  \in S^0(\R^p\times \R^r)$ and $\d_y^{2\alpha} p(y,\eta)$ is in
  $S^0(\R^p\times \R^r)$ and is thus bounded. We thus have
  \begin{align*}
    (1+|\eta|^2)^{-\hf}\ |\d_y^\alpha q (y,\eta)| 
    \leq C\: ((1+ |\eta|^2)^{-\hf}\ q (y,\eta ))^{\hf},
  \end{align*}
  which yields 
  \begin{align*}
     |\d_y^\alpha q (y,\eta)| 
    \leq C\:  (1+ |\eta|)^{\hf}\ (1+ q(y,\eta))^{\hf},
  \end{align*}
  which is the expected estimate. Let us now treat the case
  $|\beta|=1$, with for instance, $\beta=(1,0,\dots,0)$ and
  $\alpha=(0,\dots,0)$.  Define $p(y,\eta) = (1+
  |\eta|^2)^{\hf}\ q(y,\eta)$. Then $p \in S^2(\R^{p}\times
  \R^{r})$ and thus $\d_\eta^{2\beta} p(y,\eta)$ is bounded. We hence have
  \begin{align*}
    |\d_\eta^\beta p (y,\eta)|\leq C\:  (p(y,\eta))^{\hf}.
    \end{align*}
    With
    \begin{align*}
      \d_\eta^\beta p (y,\eta) = (1+
      |\eta|^2)^{\hf}\ \d_\eta^\beta q(y,\eta)  + \eta_1 (1+
      |\eta|^2)^{-\hf}\ q(y,\eta),
    \end{align*}
    the triangular inequality yields
    \begin{multline*}
      (1+
      |\eta|^2)^{\hf}\ |\d_\eta^\beta q(y,\eta)|
      \leq C\: (p(y,\eta))^{\hf} + |\eta_1| (1+
      |\eta|^2)^{-\hf}\ q(y,\eta)\\ 
      \leq C\: (q(y,\eta))^{\hf}
      ((1+
      |\eta|^2)^{\frac{1}{4}} + (q(y,\eta))^{\hf})\\
      \leq C\: (q(y,\eta))^{\hf}\ (1+ |\eta|^2)^{\frac{1}{4}}.
    \end{multline*}
    We finally obtain 
    \begin{align*}
      |\d_\eta^\beta q(y,\eta)| \leq C\: (q(y,\eta))^{\hf} (1+|\eta|)^{-\hf},
    \end{align*}
    which is the expected estimate.
\end{proof}
\begin{remark}
  \label{remark:special damping}
  If the symbol $q(z,y,\eta)$ satisfies Property (\ref{eq:P_L}) then the
  amplitude $q(z,y',\eta) +q(z,y,\eta)$ also satisfies Property (\ref{eq:P_L})
  (with derivatives \wrt $y$, $y'$ and $\eta$).
\end{remark}
\begin{proposition}
\label{Prop:exp(-Dq)}
Let $q(z,.)$ be bounded \wrt $z$ with values in $S^1(\R^p\times \R^r)$
with $q(z,.) \geq 0$. Let $q(z,.)$ satisfy Property (\ref{eq:P_L}) and
define $\rho_\Delta(z,y,\eta) = \exp[-\Delta q(z,y,\eta)]$.  Let $m\in
\N$. Then $q^m \rho_\Delta$ is smooth \wrt $\Delta$, bounded \wrt $z$,
with values in $S^0_\rho(\R^{p}\times \R^r)$ for $\Delta$ in any
interval $[\Delta_{min},\Delta_{max}]$ with $ \Delta_{min}> 0$.
\end{proposition}
\begin{proof}
  $\d_y^\alpha \d_\eta^\beta (q^m \rho_\Delta) $ is a linear combination
  of terms of the form
  \begin{align*}
    \Delta^k (\d_y^{a_1} \d_\eta^{b_1} q) \dots (\d_y^{a_l} \d_\eta^{b_l} q)
    (\d_y^{\alpha_1} \d_\eta^{\beta_1}q) \dots 
    (\d_y^{\alpha_k} \d_\eta^{\beta_k} q) 
    q^{m-l} \rho_\Delta 
  \end{align*}
  with $0\leq l\leq m$ and $a_1 + \dots +a_l + \alpha_1 + \dots
  \alpha_k = \alpha$ and $b_1 + \dots +b_l + \beta_1 + \dots \beta_k =
  \beta$. We can estimate the absolute value of each of these terms, 
  using Property (\ref{eq:P_L}), by
  \begin{multline*}
    C \Delta^k (1+ |\eta|)^{- |\beta| + \frac{|\alpha|+|\beta|}{L}}
    (1+q)^{l+k - \frac{|\alpha|+|\beta|}{L}} q^{m-l} \rho_\Delta\\
    \leq C (1+ |\eta|)^{- |\beta| + \frac{|\alpha|+|\beta|}{L}} \Delta_{min}^{-m}
  \end{multline*}
  as $(1+q)^{l+k - \frac{|\alpha|+|\beta|}{L}} q^{m-l} \Delta^{k+m}
  \rho_\Delta \leq C$.
\end{proof}

While the symbol $\exp[-\Delta q(z,y,\eta)]$ is bounded \wrt $z$
and smooth \wrt $\Delta$ with $\Delta \geq \Delta_{min}>0$ with
values in $S^0_\rho(\R^p\times\R^r)$, this fails to
be true at $\Delta=0$:
\begin{align*}
 \d_\Delta \exp[-\Delta q]|_{\Delta=0} 
 = -q \notin S^0_\rho(\R^p\times\R^r).
\end{align*}
In fact when we want to control the behavior of $\exp[-\Delta q]$ close
to $\Delta=0$ we shall use the following definition and lemmas.
\begin{definition}
  Let $L\geq 2$, $\rho =1-1/L$ and $\delta=1/L$.
  Let $\rho_\Delta(z,y,\eta)$ be a function in
  $\Cinf(\R^p\times\R^r)$ depending on the parameters
  $\Delta\geq 0$ and $z \in [0,Z]$. We say that $\rho_\Delta$
  satisfies Property (\ref{eq:Q_L}) if the following holds
  \begin{multline}
    \label{eq:Q_L}
    \tag{$Q_L$}
    \d_y^\alpha \d_\eta^\beta (\rho_\Delta -\rho_\Delta|_{\Delta=0})
    (z,y,\eta) =
    \Delta^{m+\delta(|\alpha|+|\beta|)} \rho^{m \alpha
      \beta}_\Delta(z,y,\eta), \\
    \mbox{for} \ |\alpha| + |\beta| \leq L, \ \ \
    0\leq m \leq 1 - \delta (|\alpha|+|\beta|),
  \end{multline}
  where $\rho^{m \alpha \beta}_\Delta(z,y,\eta)$ is
  bounded \wrt $\Delta$ and $z$ with values in $S^{m -
    \rho|\beta| + \delta|\alpha|}_\rho
  (\R^p\times\R^r)$.  It follows that
  $\rho_\Delta(z,y,\eta) - \rho_\Delta|_{\Delta=0}(z,y,\eta)$ is itself
  bounded \wrt $\Delta$ and $z$ with values in $S^0_\rho (\R^p\times \R^r)$.
\end{definition}
\begin{lemma}
  \label{lemma:exp(-Dq)- QL}
  Let $q(z,.)$ be bounded \wrt $z$ with values in $S^1(\R^p\times
  \R^r)$ and satisfy Property (\ref{eq:P_L}). Define
  $\rho_\Delta(z,y,\eta) = \exp[-\Delta q(z,y,\eta)]$.  Then
  $\rho_\Delta$ satisfies Property (\ref{eq:Q_L}) for $\Delta \in
  [0,\Delta_{max}]$ for any $\Delta_{max}>0$.  As
  $\rho_\Delta|_{\Delta=0} =1$, $\rho_\Delta$ is itself bounded \wrt
  $\Delta$ and $z$ with values in $S^0_\rho (\R^{p}\times \R^r)$.
\end{lemma}
\begin{proof}
  In the proof all the functions and symbols will naturally be bounded
  \wrt $z$. We thus drop the variable $z$ here for concision.

  We define 
  \begin{align*}
    \rho_\Delta^{m\alpha\beta}:= \Delta^{-m-\delta(|\alpha|+|\beta|)}
    \d_y^\alpha \d_\eta^\beta (\rho_\Delta - \rho_\Delta|_{\Delta=0}).
  \end{align*}
  We first consider the case $|\alpha|+|\beta|=0$ with $0\leq m \leq
  1$.  We need to estimate $|\d_y^a \d_\eta^b \rho_\Delta^{m00}|$.  The
  case $m=0$, $|a+b|=0$ has to be treated independently but is
  trivial: we clearly have $|\rho_\Delta^{000}|= |\rho_\Delta -1| \leq
  C$.  We shall now estimate $|\d_y^a \d_\eta^b \rho_\Delta^{m00}| = |
  \Delta^{-m} \d_y^a \d_\eta^b ( \rho_\Delta-1 )|$ in the case where
  $m>0$ or $|a|+|b|>0$.  For this we write 
  \begin{align}
    \label{eq:Taylor rho_Delta}
    \rho_\Delta(y,\eta) -1 = -\Delta \int_0^1 q(y,\eta)
    \exp[-s\Delta q(y,\eta)] d s.
  \end{align}
  We then have
  $\rho_\Delta^{m00}(y,\eta) = -\int_0^1 q_\Delta^m(s,y,\eta) ds$ with
  \begin{align*}
  q_\Delta^m(s,y,\eta) = \Delta^{1-m} q(y,\eta) \exp[-s\Delta q(y,\eta)].
  \end{align*}
  We prove that
  \begin{align*}
    |\d_y^a \d_\eta^b q_\Delta^m(s,y,\eta)| \leq C(s)
    (1+|\eta|)^{m-\rho |b|+ \delta |a|}
  \end{align*}
  with $C(s)$ bounded \wrt $\Delta$ and $L^1$ \wrt $s \in [0,1]$.
  The result then follows for $\rho_\Delta^{m00}$.

  When computing $\d_y^a\d_\eta^b q_\Delta^m$ we obtain
  a linear combination of terms of the form
  \begin{multline*}
    \Delta^{1-m} (\d^{a_0}_y\d^{b_0}_\eta q)
    (-s \Delta)^k (\d^{a_1}_y\d^{b_1}_\eta q)\dots
    (\d^{a_k}_y\d^{b_k}_\eta q)
    \exp[-s\Delta q],\\
    \mbox{with}\ a_0 + a_1 + \dots + a_k = a, \ \ 
    b_0 + b_1 + \dots + b_k = b,
  \end{multline*}
  (where $k$ can be $0$). Using Property (\ref{eq:P_L}) we find that
  the absolute value of such a term is bounded by
  \begin{multline*}
    C \Delta^{1-m}(s \Delta)^k (1+|\eta|)^{-|b| + \delta(|a+b|)}
    (1+q)^{k+1-\delta(|a+b|)} \exp[-s\Delta q]\\
    \leq C s^{m +\delta(|a+b|) -1}(1+|\eta|)^{m -\rho|b|+\delta|a|}
    \Delta^{\delta(|a+b|)} \\
    (s \Delta (1+q))^{-m + k+1-\delta(|a+b|)}\exp[-s\Delta q],
  \end{multline*}
  as $1\leq C(1+|\eta|)^m (1+q)^{-m}$ if $m\geq 0$.  If $l:=-m +
  k+1-\delta(|a+b|) \geq 0$ we use that $(s\Delta (1+q))^l
  \exp[-s\Delta q]\leq C$ if $0\leq s \leq 1$, $0 \leq \Delta \leq
  \Delta_{max}$ and $q\geq 0$ and we obtain the following estimate
  \begin{align*}
    C s^{m +\delta(|a+b|) -1}(1+|\eta|)^{m -\rho|b| + \delta|a|}
    \Delta^{\delta(|a+b|)}.
  \end{align*}
  If $l<0$, $(1+q)^l$ is simply bounded ($q\geq 0$) and we obtain the
  following estimate:
  \begin{align*}
    C \Delta^{k+1-m} s^k (1+|\eta|)^{m -\rho|b| + \delta |a|}.
  \end{align*}
  As $m +\delta(|a+b|) -1>-1$ in the considered case, both estimates
  exhibit bounds that are in $L^1([0,1])$ \wrt $s$. We also
  have  uniform bounds \wrt $\Delta$ as we have assumed $m\leq 1$.
  

  We now treat the case $1\leq |\alpha| + |\beta| \leq L$, $0\leq m
  \leq 1 - \delta (|\alpha|+|\beta|)$.  We estimate the absolute value
  of
  \begin{align*}
    \d_y^a \d_\eta^b( \rho_\Delta^{m\alpha\beta}) =
    \Delta^{-m-\delta(|\alpha|+|\beta|)}\
    \d_y^{a+\alpha} \d_\eta^{b+\beta}\rho_\Delta 
  \end{align*}
  which is a linear combination of terms of the form
  \begin{multline*}
    \Delta^{k-m-\delta(|\alpha|+|\beta|)}\ (\d^{a_1}_y\d^{b_1}_\eta
    q)\dots (\d^{a_k}_y\d^{b_k}_\eta q)
    \exp[-\Delta q],\\
    \mbox{with}\ a_1 + \dots + a_k = a+\alpha, \ \ b_1 + \dots + b_k =
    b+\beta,
  \end{multline*}
  where $k\geq 1$. Using Property (\ref{eq:P_L}) we
  find that the absolute value of such a term is bounded by
  \begin{multline*}
    C\Delta^{k-m-\delta(|\alpha|+|\beta|)}
    (1+|\eta|)^{-|\beta|-|b| + \delta(|\alpha|+|a|+|\beta|+|b|)} \\
    (1+q)^{k -\delta(|\alpha|+|a|+|\beta|+|b|)} \exp[-\Delta q]\\
    \leq C (1+|\eta|)^{m - \rho (|\beta|+|b| ) + \delta (|\alpha|+|a|)}
    (1+q)^{-\delta (|a|+|b|)}
    \\
    (\Delta (1+q))^{k -m - \delta(|\alpha|+|\beta|)}\exp[-\Delta q]\\
    \leq C (1+|\eta|)^{m - \rho (|\beta|+|b| ) + \delta (|\alpha|+|a|)},
  \end{multline*}
  as $k -m - \delta (|\alpha|+|\beta|) \geq 1- m -
  \delta(|\alpha|+|\beta|)\geq 0$ and $0\leq \Delta \leq
  \Delta_{max}$. This completes the proof.
\end{proof}
\begin{lemma}
  \label{lemma:QL stable by f}
  Let $f \in \Cinf(\R)$ and $q_\Delta(z,y,\eta)$ in $\Cinf(\R^p\times
  \R^r)$ that satisfies Property~(\ref{eq:Q_L}) and such that
  $q_\Delta(z,.)|_{\Delta=0}$ is independent of $y$ and $\eta$. Then
  $f(q_\Delta)(z,y,\eta)$ satisfies Property (\ref{eq:Q_L}).
\end{lemma}
\begin{proof}
  Again bounds \wrt $z$ are clear.
  We first treat the case $|\alpha|+|\beta|=0$. We write
  \begin{align*}
    f(q_\Delta) - f(q_\Delta|_{\Delta=0}) = (q_\Delta -q_\Delta|_{\Delta=0})
    \int_0^1 f'((1-s)q_\Delta|_{\Delta=0} + s q_\Delta) d s.
  \end{align*}
  As $q_\Delta|_{\Delta=0}$ is independent of $y$ and $\eta$, then
  $q_\Delta$ is bounded \wrt $\Delta$ with values in
  $S^0_\rho(\R^p\times \R^r)$ by Property~(\ref{eq:Q_L})  and so are
  $(1-s)q_\Delta|_{\Delta=0} + s q_\Delta$ and
  $f'((1-s)q_\Delta|_{\Delta=0} + s q_\Delta)$ by Lemma 18.1.10 in
  \cite{Hoermander:V3} with bounds in $S^0_\rho(\R^p\times
  \R^r)$ uniform with respect to $s$. We thus obtain that $\int_0^1
  f'((1-s)q_\Delta|_{\Delta=0} + s q_\Delta) d s$ is bounded \wrt
  $\Delta$ with values in $S^0_\rho(\R^p\times \R^r)$. We
  conclude using Property (\ref{eq:Q_L}) for $q_\Delta
  -q_\Delta|_{\Delta=0}$.  Let us now treat the case $1 \leq
  |\alpha|+|\beta|\leq L$ and choose $0\leq m \leq 1 -
  \delta(|\alpha|+|\beta|)$.  We see that $\d_y^\alpha \d_\eta^\beta
  f(q_\Delta) $ is a linear combination of terms of the form
  \begin{align*}
    (\d_y^{\alpha_1} \d_\eta^{\beta_1} q_\Delta)\dots 
    (\d_y^{\alpha_k} \d_\eta^{\beta_k} q_\Delta) f^{(k)}(q_\Delta),
  \end{align*} 
  where $k \geq 1$, $\alpha_1 + \dots + \alpha_k = \alpha$, $\beta_1+
  \dots + \beta_k = \beta$.  Now choose $0\leq m_i \leq 1 -
  \delta(|\alpha_i|+|\beta_i|)$, $i=1,\dots, k$, such that $m=m_1+
  \dots +m_k$. Then Property (\ref{eq:Q_L}) yields terms of the form
  \begin{align*}
    \Delta^{m_1 + \delta(|\alpha_1|+|\beta_1|)}\dots \Delta^{m_k +
      \delta(|\alpha_k|+|\beta_k|)} q_\Delta^{m_1 \alpha_1 \beta_1 }
    \dots q_\Delta^{m_k \alpha_k \beta_k} = \Delta^{m + \delta(|\alpha|+|\beta|)}
    q_\Delta^{m\alpha \beta}
  \end{align*}
  with $q_\Delta^{m_i \alpha_i \beta_i}$, $i=1,\dots, k$, bounded
  \wrt $\Delta$ with values in $S^{m_i - \rho|\alpha_i| + \delta
    |\beta_i| }_\rho(\R^p\times \R^r)$ and $q_\Delta^{m\alpha \beta} :=
  q_\Delta^{m_1 \alpha_1 \beta_1 } \dots q_\Delta^{m_k \alpha_k
    \beta_k}$. We note that $f^{(k)}(q_\Delta)$ is bounded \wrt
  $\Delta$ with values in $S^0_\rho(\R^p\times \R^r)$.  The symbol
  $q_\Delta^{m\alpha \beta}$ is bounded \wrt  $\Delta$ with values
  in $S^{m - \rho|\alpha| + \delta |\beta| }_\rho(\R^p\times \R^r)$, which
  yields the result.
\end{proof}
With Remark~\ref{remark:special damping}, Lemma~\ref{lemma:QL stable
by f} and the previous lemma we obtain
\begin{corollary}
  \label{corollary: property QL for p_Delta}
  Let $f \in \Cinf(\R)$ and let $q(z,.)$ bounded \wrt $z$
  with values in $S^1(\R^{p}\times \R^r)$ satisfy Property
  (\ref{eq:P_L}). Define
  \begin{align*}
  p_\Delta(z,y',y,\eta) =\exp[-\Delta (q(z,y',\eta)+ q(z,y,\eta))].
  \end{align*}
  Then $f(p_\Delta)$ satisfies Property (\ref{eq:Q_L}). As
  $f(p_\Delta)|_{\Delta=0}=1$, $f(p_\Delta)$ is itself bounded \wrt
  $\Delta$ and $z$ with values in $S^0_\rho(\R^{2p}\times
  \R^r)$.
\end{corollary}
Note that the property (\ref{eq:Q_L}) is stable when we go from
amplitudes to symbols:
\begin{proposition}
  Let $q_\Delta(z,x,y,\xi)$ be an amplitude in
  $S^0_\rho(\R^{2p}\times \R^p)$ depending on the
  parameters $\Delta\geq 0$ and $z \in [0,Z]$ that satisfies Property
  (\ref{eq:Q_L}). Then $\symb{q_\Delta} (z,x,\xi)$ satisfies property
  (\ref{eq:Q_L}).
\end{proposition}
\begin{proof}
  We use the oscillatory integral representation for the symbol:
  \begin{align*}
    \symb{q_\Delta}(z,x,\xi) := \iint \exp[-i \inp{y}{\eta}]\
    q_\Delta(z,x,x-y,\xi-\eta) \ \dslash \eta \ d y.
  \end{align*}
  Let $0\leq |\alpha|+|\beta|\leq L$ and $0\leq m \leq 1 - \delta
  (|\alpha|+|\beta|)$.  Computing $\d_x^\alpha \d_\xi^\beta (
  \symb{q_\Delta} - \symb{q_\Delta}|_{\Delta=0})$, we obtain a linear
  combination of terms of the form, with $\alpha_1+\alpha_2 = \alpha$,
  \begin{multline*}
    \iint \exp[-i \inp{y}{\eta}]\ \d_2^{\alpha_1} \d_3^{\alpha_2}
    \d_4^\beta (q_\Delta-q_\Delta|_{\Delta=0}) (z,x,x-y,\xi-\eta)\
    \dslash \eta \ d y\\ = \iint \exp[-i \inp{y}{\eta}]\
    \Delta^{m+\delta(|\alpha|+|\beta|)} q_\Delta^{m
    (\alpha_1,\alpha_2) \beta}(z,x,x-y,\xi-\eta) \ \dslash \eta \ d
    y\\ = \Delta^{m+\delta(|\alpha|+|\beta|)} \symb{q_\Delta^{m
    (\alpha_1,\alpha_2) \beta}},
  \end{multline*}
  where $q_\Delta^{m (\alpha_1,\alpha_2) \beta}$ is bounded \wrt
    $\Delta$ and $z$ with values in the symbol space $S^{m -
    \rho|\beta| + \delta|\alpha|}_\rho (\R^{2p}\times \R^p)$.  As the
    map $a \mapsto \symb{a}$ maps bounded sets into bounded sets the
    result follows.
\end{proof}
We shall also need the following lemma.
\begin{lemma}
  \label{lemma:symb - amp}
  Let $q_\Delta(z,x,y,\xi)$ be an amplitude in
  $S^0_\rho(\R^{2p}\times \R^p)$ depending on
  the parameters $\Delta\geq 0$ and $z \in [0,Z]$ that satisfies
  Property (\ref{eq:Q_L}) for $1\leq |\alpha|+|\beta|\leq 2$ and such that
  $q_\Delta(z,.)|_{\Delta=0}$ is independent of $(x,y,\xi)$. Let
  $r(x,\xi) \in S^s(\R^p\times \R^p)$ for some $s \in \R$.  Then
  \begin{align*}
    \symb{q_\Delta\ r } (z,x,\xi) - q_\Delta(z,x,x,\xi)\ r(x,\xi) =
    \Delta^{m + 2 \delta} \lambda_{\Delta }^{m}(z,x,\xi), \ 0\leq m \leq \rho
    - \delta,
  \end{align*}
  where the function $\lambda_\Delta^{m }(z,x,\xi)$ is bounded with
  respect to $\Delta$ and $z$ with values in $S^{m +s - (\rho -
  \delta)}_\rho(\R^p\times \R^p)$.
\end{lemma}
\begin{proof}
  For the sake of concision we take $p=1$ in the proof but it
  naturally extends to $p\geq 1$. We write $\lambda_\Delta = q_\Delta
  r$. Using the oscillatory integral representation of
  $\symb{q_\Delta}$ we obtain
  \begin{multline*}
    \symb{q_\Delta r } (z,x,\xi) - q_\Delta(z,x,x,\xi) r(x,\xi) \\
    =\iint \exp[-i \inp{y}{\xi-\eta}]
    (\lambda_\Delta(z,x,x-y,\eta)  - \lambda_\Delta(z,x,x,\eta)\
    \dslash \eta \ d y.
  \end{multline*}
  Taylor's formula yields
  \begin{multline*}
    \symb{q_\Delta r } (z,x,\xi) - q_\Delta(z,x,x,\xi)r(x,\xi) \\
    =\int_0^1 \iint -y \exp[-i \inp{y}{\xi-\eta}] \
   \d_3 \lambda_\Delta(z,x,x-s y,\eta) \
    \dslash \eta \ d y\ d s.
  \end{multline*}
  With an integration by parts we obtain
  \begin{multline*}
    \symb{q_\Delta r} (z,x,\xi) - q_\Delta(z,x,x,\xi)r(x,\xi) \\
    =\int_0^1 \iint i  \exp[-i \inp{y}{\xi-\eta}] \ \d_3 \d_4
    \lambda_\Delta(z,x,x-s y,\eta) \ \dslash \eta \ d y\ d s\\
    = \symb{i \int_0^1 \d_3 \d_4 \lambda_\Delta(z,x,(1-s)x+s y,\xi)\ ds},
  \end{multline*}
  where $\d_3 \d_4 \lambda_\Delta(z,x,y,\xi) = (\d_y \d_\xi
  q_\Delta)(z,x,y,\xi) \ r(x,\xi) + \d_y q_\Delta(z,x,y,\xi) \d_\xi
  r(x,\xi)$, as $r$ does not  depend on $y$.  The first term is treated
  using Property (\ref{eq:Q_L}) while for the second one we write
  \begin{align*}
    \d_y q_\Delta\  \d_\xi r 
    = \Delta^{m'+\delta}\  q_\Delta^{m'(0,1)0}\  \d_\xi r,
  \end{align*}
  where $0\leq m' \leq 1-\delta$ and $q_\Delta^{m'(0,1)0}\ r \in
  S^{m'+s-1+\delta}_\rho(\R^{2p}\times \R^p)$ by Property
  (\ref{eq:Q_L}). We actually take $\delta \leq m' \leq 1-\delta$ and
  write $m=m'-\delta$.  We obtain
  \begin{align*}
    \d_y q_\Delta \d_\xi r = \Delta^{m+2\delta}\  \tilde{q}_\Delta^{m},
  \end{align*}
  where $\tilde{q}_\Delta^{m}$ is bounded \wrt $\Delta$ with values
  in $S^{m+s-\rho+\delta}_\rho(\R^{2p}\times \R^p)$ and
  $0\leq m\leq 1-2\delta= \rho -\delta$.  We conclude since the map
  $\sigma\{.\}$ maps bounded sets into bounded sets.
\end{proof}
We are now ready to give an estimate of the operator norm of the
thin-slab propagator, $\G_{(z',z)}$, in $L(H^{(s)}(X),H^{(s)}(X'))$
for any $s \in \R$.
\addtocounter{delta}{1}
\begin{theorem}
  \label{theorem:H^s estimate}
  Let $s \in \R$. There exists $M>0$, $\Delta_\thedelta
  >0$ such that
  \begin{align*}
  \| \G_{(z',z)}\|_{(H^{(s)},H^{(s)})}\leq 1+ \Delta M,
  \end{align*}
  for all $z',z \in [0,Z]$ such that $0 \leq\Delta=z'-z\leq
  \Delta_\thedelta$.
\end{theorem}
In the proof we assume that $c_1$ satisfies property $P_L$ for some $L
\geq 2$. We know that it is always true for $L=2$ by Lemma~\ref{lemma:
PL with L=2 for c_1} but special choices for $c_1$ can be made. As
before we use $\rho=1-1/L$ and $\delta=1/L$ with $\rho > \delta$ for
$L>2$ and $\rho=\delta=\hf$ for $L=2$. In the proof we proceed
classically by computing $\G_{(z',z)} \G_{(z',z)}^\ast$ and use the
classical results on $\psi$DOs (see e.g. \cite[Section 5]{MeSj:76} and
also \cite{Hoermander:83}). Here we however do not content ourself
with the continuity of $\G_{(z',z)}$ but we want to obtain a precise
estimate of the operator norm in $L(H^{(s)}(X),H^{(s)}(X'))$, which
will be required in the sequel. Here we exploit the fact that
$\Delta$ can be taken arbitrarily small  which allows to carry out
some explicit computations.
\begin{proof}
  Let $s \in \R$, then the kernel of $\A_{(z',z)}:=\G_{(z',z)}\circ
  E^{(-s)}$ is given by
  \begin{align*}
    \A_{(z',z)}(x',x)
    = \int \exp[i
    \phi_{(z',z)}(x',x,\xi)]\ g_{(z',z)}(x',\xi)\
    \e{\xi}{-s}\, \dslash \xi.
  \end{align*}
  Computing the kernel $D_{(z',z)}$ of ${\cal
    D}_{(z',z)}:=\A_{(z',z)}\circ\A_{(z',z)}^\ast$ we obtain
  \begin{multline*}
    D_{(z',z)}(x',x) \\
    = \int \exp\left[ i\inp{x'-x}{\xi}
    + i\Delta \left(b_1(z,x',\xi)-b_1(z,x,\xi)\right)\right]
    d_{(z',z)}(x',x,\xi)\
    \dslash \xi
  \end{multline*}
  where
  \begin{multline*}
    d_{(z',z)}(x',x,\xi) \\
    = \exp[-\Delta (c_1(z,x',\xi)+c_1(z,x,\xi) ] g_{(z',z)}(x',\xi)\ 
    \overline{g_{(z',z)}}(x,\xi)\ \e{\xi}{-2s}.
  \end{multline*}
  We write $b_1(z,x',\xi)-b_1(z,x,\xi) = \inp{x'-x}{ h(z,x',x,\xi)}$
  where $h$ is smooth and homogeneous of degree one in $\xi$,
  $|\xi|\geq 1$.  The function $h$ and continuous \wrt $z$ with
  values in $S^1(X\times \R^n)$ by Assumption~\ref{assumpt:general
    assumption} and estimate (1.1.9) in \cite{Hoermander:V1}. We thus
  obtain that the change of variables $\xi \to \xi+\Delta
  h(z,x',x,\xi)=H_{(\Delta,z,x',x)}(\xi)$ is a global diffeomorphism
  for $\Delta$ small enough (uniformly in $z \in [0,Z]$).  We denote
  $\tilde{\xi} (\Delta,z,x',x,\xi) = H^{-1}_{(\Delta,z,x',x)}(\xi)$.
  We thus have
  \begin{align*}
    D_{(z',z)}(x',x) =\!
    \int \exp\left[ i\inp{x'-x}{\xi} \right]
    d_{(z',z)}(x',x,\tilde{\xi}(\Delta,z,x',x,\xi))\
    \J_\Delta(z,x',x,\xi)\  \dslash \xi
  \end{align*}
  where $\J_\Delta(z,x',x,\xi)$ is the Jacobian.
  \begin{lemma}
    \label{lemma: change of variables}
    The function $\tilde{\xi}(\Delta,z,x',x,\xi)$ is homogeneous of degree 1 in
    $\xi$, for $|\xi|\geq 1$,  continuous \wrt $z$, $\Cinf$ \wrt $\Delta$ with values in
    $S^1(\R^{2n}\times(\R^{n}))$ if
    $\Delta$ is small enough, i.e.,
    \begin{align*}
    \exists \Delta_\thedelta >0, \ \
    \tilde{\xi} \in \Con^0([0,Z],\Cinf([0,\Delta_\thedelta],
    S^1(\R^{2n}\times(\R^{n})))).
    \end{align*}
  \end{lemma}
  \begin{proof}
    Homogeneity is clear. We
    have
    \begin{multline*}
      |\tilde{\xi}(\Delta,z,x',x,\xi)| = |\xi - 
      \Delta h (z,x',x,\tilde{\xi}(\Delta,z,x',x,\xi))| \\
      \leq 1 + \Delta C (1+|\tilde{\xi}(\Delta,z,x',x,\xi)|),\ \ |\xi|=1,
    \end{multline*}
    which yields, because of homogeneity,
    \begin{align*}
      |\tilde{\xi}(\Delta,z,x',x,\xi)| \leq \frac{1+ \Delta C}{1-
        \Delta C}(1+ |\xi|), \ |\xi|\geq 1, 
    \end{align*}
    for $\Delta$ small enough, uniformly chosen \wrt $z \in
    [0,Z]$, $x',x \in \R^n$.  Differentiating the $j^{th}$ coordinate
    of $\xi$,
    \begin{align*}
    \xi_j = \tilde{\xi}_j(\Delta,z,x',x,\xi) + \Delta
    h_j (z,x',x,\tilde{\xi}(\Delta,z,x',x,\xi)),\ j=1,\dots,n,
    \end{align*}
    \wrt $x_i$ yields
    \begin{multline}
      \label{eq:dx eta}
      \d_{x_i} \tilde{\xi}_j(\Delta,z,x',x,\xi) +
      \Delta \partial_{x_i} h_j (z,x',x,\tilde{\xi}(\Delta,z,x',x,\xi))\\
      +\Delta \sum_l \partial_{\tilde{\xi}_l}
      h_j (z,x',x,\tilde{\xi}(\Delta,z,x',x,\xi)) \ 
      \d_{x_i} \tilde{\xi}_l(\Delta,z,x',x,\xi) = 0, \\ j=1,\dots,n.
    \end{multline}
    The partial derivatives of $h$ are bounded for $|\xi|=1$.
    We can hence solve for $\d_{x_i} \tilde{\xi}(\Delta,z,x',x,\xi)$
    when $\Delta$ is small enough and find the expected estimate from
    that obtained for $\tilde{\xi}(\Delta,z,x',x,\xi)$:
    \begin{align*}
    \exists C> 0, \  |\d_{x_i} \tilde{\xi}(\Delta,z,x',x,\xi)|
    \leq C (1+ | \xi|),
    \ x',x\in\R^{n}, \ \xi \in \R^n.
    \end{align*}
    Differentiating \wrt $x_i'$, $\xi_i$, and
    $\Delta$ yields similar structures and the proper symbol
    estimates. The proof carries on by induction. Note that the
    required size for $\Delta$ to solve the systems of the form
    (\ref{eq:dx eta}) remains fixed along the induction process.
  \end{proof}
  {\em Continuation of the proof of Theorem~\ref{theorem:H^s
      estimate}.}  From (the proof of) Lemma~\ref{lemma: change of
      variables} we also obtain that the Jacobian
      $\J_\Delta(z,x',x,\xi)$ is homogeneous of degree zero in $\xi$,
      $|\xi|\geq 1$, and is continuous \wrt $z$ and $\Cinf$ \wrt
      $\Delta$ with values in $S^0(\R^{2n}\times \R^{n}))$.
  
  We write $\tilde{p}_\Delta(z,x',x,\xi):= \exp[-\Delta
  (c_1(z,x',\xi)+c_1(z,x,\xi) ]$.  As $c_1$ satisfies Property
  (\ref{eq:P_L}) we then have $\tilde{p}_\Delta$ satisfying property
  (\ref{eq:Q_L}) by Corollary~\ref{corollary: property QL for
  p_Delta}.  Define $p_\Delta(z,x',x,\xi) :=
  \tilde{p}_\Delta(z,x',x,\tilde{\xi}(\Delta,z,x',x,\xi))$.
  Lemma~\ref{lemma: change of variable - symbol} and Lemma~\ref{lemma:
  change of variables} yield that $p_\Delta$ satisfy property
  (\ref{eq:Q_L}). We then have
  \begin{align*}
    d_{(z',z)}(x',x,\tilde{\xi}(\Delta,z,x',x,\xi))\
    \J_\Delta(z,x',x,\xi)\
    =:  p_\Delta(z,x',x,\xi)\ k_\Delta(z,x',x,\xi)
  \end{align*}
  where $k_\Delta(z,.)$ is bounded \wrt $z$ and $\Cinf$
  \wrt $\Delta$ with values in $S^{-2s} (\R^{2n} \times
  \R^n)$ and $k_\Delta(z,.)|_{\Delta=0}= \e{.}{-2s}$ by
  Lemma~\ref{lemma: change of variable - symbol} and Lemma~\ref{lemma:
    change of variables}. By Theorem 1.1.9 and formula (1.1.9) in
  \cite{Hoermander:V1} we obtain
  \begin{align*}
    k_\Delta(z,x',x,\xi) = \e{\xi}{-2s} + \Delta
    \tilde{k}_\Delta(z,x',x,\xi),
  \end{align*}
  where $\tilde{k}_\Delta$ is bounded \wrt $z$ and $\Cinf$
  \wrt $\Delta$ with values in $S^{-2s} (\R^{2n} \times
  \R^n)$.

  Call ${\cal F}_{(z',z)}:= E^{(s)} \circ {\cal D}_{(z',z)}\circ
  E^{(s)}$.  Its symbol is in $S^0_\rho(\R^{n} \times \R^n)$
  and is given by
  \begin{multline*}
    f_{(z',z)}(x',\xi):= (\e{\xi}{s}\ \#\
    \symb{ p_\Delta(z,x',x,\xi) \
      k_\Delta(z,x',x,\xi) }\
    \#\ \e{\xi}{s}) (x',\xi)\\
    = (\e{\xi}{s}\ \#\ \symb{ p_\Delta(z,x',x,\xi)
      \e{\xi}{-2s}}\ \#\ \e{\xi}{s}) (x',\xi)\\
    + \Delta (\e{\xi}{s}\ \#\ \symb{ p_\Delta(z,x',x,\xi)
       \tilde{k}_\Delta(z,x',x,\xi) \e{\xi}{-2s} }\ \#\ \e{\xi}{s} ) (x',\xi)
  \end{multline*}
  As $p_\Delta$ bounded \wrt $z$ and $\Delta$, $\Delta$
  small enough, with values in $S^0_\rho(\R^{2n} \times \R^n)$
  (Property (\ref{eq:Q_L})) we obtain that the second term in the equation
  above satisfies the same property and thus we can write
  \begin{align*}
    {\cal F}_{(z',z)}={\cal F}^a_{(z',z)}+ \Delta {\cal F}^1_{(z',z)}
  \end{align*}
  where ${\cal F}^a_{(z',z)}$ has for symbol
  \begin{align*}
    (\e{\xi}{s}\  \#\ \symb{ p_\Delta(z,x',x,\xi)
      \e{\xi}{-2s} }\ \#\ \e{\xi}{s} ) (x',\xi)
  \end{align*}
  and $\|{\cal F}^1_{(z',z)}\|_{(L^2,L^2)} \leq K^1$, uniformly in $z
  \in [0,Z]$ and $\Delta$, $\Delta$ small enough, by the
  Calder\'{o}n-Vaillancourt theorem (see \cite[Chapter 7, Sections
  1,2]{Kumano-go:81} or \cite[Section XIII-2]{Taylor:81}) in the case
  $L=2$ and by Theorem 18.1.11 in \cite{Hoermander:V3} in the case
  $L>2$.  With Lemma~\ref{lemma:symb - amp} we see that
  \begin{align*}
    \symb{ p_\Delta(z,x',x,\xi)
     \e{\xi}{-2s} } - p_\Delta(z,x',x',\xi)
    \e{\xi}{-2s} = \Delta \lambda_\Delta (z,x',\xi)
  \end{align*}
  where $\lambda_\Delta$ is bounded \wrt $\Delta$ and $z$
  with values in $S^{-2s}_\rho(\R^{n} \times \R^n)$.
  We thus obtain
  \begin{align*}
    {\cal F}^a_{(z',z)}={\cal F}^b_{(z',z)}+ \Delta {\cal F}^2_{(z',z)}
  \end{align*}
  where ${\cal F}^b_{(z',z)}$ has for symbol
  \begin{multline*}
    f_\Delta^b(z,x',\xi) :=
    (\e{\xi}{s}\ \#\ p_\Delta(z,x',x',\xi) \e{\xi}{-2s}\ \#\ \e{\xi}{s})
    (z,x',\xi) \\
    = (\e{\xi}{s}\ \#\ p_\Delta(z,x',x',\xi)  \e{\xi}{-s})(z,x',\xi)
  \end{multline*}
  and $\|{\cal F}^2_{(z',z)}\|_{(L^2,L^2)} \leq K^2$ uniformly in $z
  \in [0,Z]$ and $\Delta$, $\Delta$ small enough.

  For the rest of the proof, if we don't write it explicitly, by
  $p_\Delta$ and $p_\Delta(z,x,\xi)$ we shall actually mean
  $p_\Delta(z,x,x,\xi)$.
  \begin{lemma}
    \begin{align*}
      (\e{.}{s}\ \#\ p_\Delta(z,.)  \e{.}{-s})(z,x,\xi)-
      p_\Delta(z,x,\xi) = \Delta \mu_\Delta(z,x,\xi),
    \end{align*}
    where $\mu_\Delta(z,x,\xi)$ is bounded \wrt $z$
    and $\Delta$ with values in $S^{0}_\rho(X\times \R^n)$.
  \end{lemma}
  \begin{proof}
    We write
    \begin{align*}
      p_\Delta(z,x,\xi) = \e{\xi}{-s}
      \iint\exp[-i\inp{y}{\xi-\eta}]\
      \e{\eta}{s} \ p_\Delta(z,x,\xi)\ \dslash
      \eta\ d y
    \end{align*}
    and thus obtain, with the oscillatory integral representation for
    the composition formula,
    \begin{multline*}(\e{.}{s}\ \#\ p_\Delta(z,.)  \e{.}{-s})(z,x,\xi)-
      p_\Delta(z,x,\xi) =\\
      \e{\xi}{-s}
      \iint\exp[-i\inp{y}{\xi-\eta}]\
      \e{\eta}{s}\
      (p_\Delta(z,x-y,\xi) - p_\Delta(z,x,\xi))\
      \dslash \eta \ d y.
    \end{multline*}
    With Taylor's formula and applying an integration by part, we find
    (we have supposed $n=1$ for the sake of simplicity but it
    naturally extends to $p\geq1$)
    \begin{multline*}(\e{.}{s}\ \#\ p_\Delta(z,.)  \e{.}{-s})(z,x,\xi)-
      p_\Delta(z,x,\xi) =\\
      \e{\xi}{-s}
      \int_0^1 \iint i \exp[-i\inp{y}{\xi-\eta}]\
      \d_\eta \e{\eta}{s}\
      \d_x p_\Delta(z,x- r y,\xi) \
      \dslash \eta\ d y\ d r .
    \end{multline*}
    Using Property (\ref{eq:Q_L}) with $m=1-\delta$ we find
    \begin{multline*}(\e{.}{s}\ \#\ p_\Delta(z,.)  \e{.}{-s})(z,x,\xi)-
      p_\Delta(z,x,\xi) = \Delta \e{\xi}{-s}\\
      \int_0^1 \iint i \exp[-i\inp{y}{\xi-\eta}]\
      \d_\eta \e{\eta}{s}\
      q_\Delta^{m 1 0}(z,(1-r)x+r (x-y),\xi) \
      \dslash \eta\ d y\ d r \\
      = \Delta \e{\xi}{-s}
      (\d_\xi \e{\xi}{s}\ \#\ \tilde{q}_{\Delta}^{m 1 0}(z,u,x,\xi))|_{u=x}
    \end{multline*}
    where
    \begin{align*}
      \tilde{q}_{\Delta}^{m 1 0}(z,u,x,\xi) =
      \int_0^1 q_\Delta^{m 1 0}(z,(1-r)u+ r x,\xi) \ d r.
    \end{align*}
    As $\tilde{q}_{\Delta}^{m 1 0}$ is bounded \wrt
    $\Delta$ and $z$ with values in $S^1_\rho(\R^n\times \R^n)$
    we obtain the result.
  \end{proof}
  {\em End of the proof of Theorem~\ref{theorem:H^s
      estimate}.}
  With the previous lemma we see that
  \begin{align*}
    {\cal F}^b_{(z',z)}={\cal F}^c_{(z',z)}+ \Delta {\cal F}^3_{(z',z)}
  \end{align*}
  where ${\cal F}^c_{(z',z)}$ has for symbol
  $p_\Delta(z,x,x',\xi)$ and $\|{\cal
    F}^3_{(z',z)}\|_{(L^2,L^2)} \leq K^3$ uniformly in $z \in [0,Z]$
  and $\Delta$, $\Delta$ small enough.

  To estimate $\|{\cal F}^c_{(z',z)}\|_{(L^2,L^2)}$ we follow the
  procedure at the end of the proof of Theorem 18.1.11 in
  \cite{Hoermander:V3}. Let $A:=1+\Delta$. Define
  \begin{align*}
    \nu_\Delta(z,x',\xi) =
    \sqrt{A^2 - |p_\Delta(z,x',x',\xi)|^2},
  \end{align*}
  which satisfies Property (\ref{eq:Q_L}) by Lemma~\ref{lemma:QL stable by
    f}. Then define $r_\Delta$ by
  \begin{align*}
    \nu_\Delta\ \#\ \nu_\Delta^\ast = A^2 - p_\Delta\ \#\ p_\Delta^\ast - r_\Delta.
  \end{align*}
  Note that
  \begin{align*}
    \nu_\Delta\ \#\ \nu_\Delta^\ast(z,x',\xi) =
    \symb{\nu_\Delta (z,x',\xi) \overline{\nu_\Delta}(z,x,\xi)}
    (z,x,\xi)
  \end{align*}
  It is easy to check that $\nu_\Delta (z,x',\xi)
  \overline{\nu_\Delta}(z,x,\xi)$ satisfies Property (\ref{eq:Q_L})
  for $|\alpha| + |\beta|\geq 1$. The same applies to $p_\Delta
  (z,x',\xi) \overline{p_\Delta}(z,x,\xi)$.  Lemma~\ref{lemma:symb -
  amp} applies and with $m=\rho-\delta$ we thus obtain that $r_\Delta
  = \Delta \tilde{r}_\Delta$ with $\tilde{r}_\Delta$ bounded \wrt $z$
  and $\Delta$ with values $S^0_\rho(X\times \R^n)$.  Thus
  \begin{align*}
    \|{\cal F}^c_{(z',z)}\|_{(L^2,L^2)} =
    \|({\cal F}^c_{(z',z)})^\ast\|_{(L^2,L^2)}
    \leq \sqrt{(1+\Delta)^2 + \Delta C}\leq 1+ \Delta K^4,
  \end{align*}
  for some $K^4>0$ large enough.  We thus obtain that $\|{\cal
    F}_{(z',z)}\|_{(L^2,L^2)} \leq 1 + \Delta K$ where $K=K^1 + K^2 +
  K^3 + K^4$. With the definition of ${\cal F}_{(z',z)}$ it follows
  that
  \begin{align*}
    \|\G_{(z',z)}\|_{(H^{(s)},H^{(s)})} =
    \|(\G_{(z',z)})^\ast\|_{(H^{(s)},H^{(s)})} \leq \sqrt{1 + \Delta K}
  \end{align*}
  which concludes the proof of Theorem~\ref{theorem:H^s estimate}.
\end{proof}

We observe that for $\Delta$ small enough, the function $\inp{x'}{\xi}
+ \Delta b_1(z,x',\xi)$ satisfies the conditions $(P)$-$(i)$,
$(P)$-$(ii)$, and $(P)$-$(iii)$ in \cite[page 2]{Kumano-go:76}.  With
Lemmas~\ref{lemma: PL with L=2 for c_1} and ~\ref{lemma:exp(-Dq)- QL},
we observe that an FIO with phase function $\phi_{(z',z)}(x',x,\xi)$
and amplitude in $\sigma_A(z,x',\xi)$ in $S^m(X\times \R)$ may
actually be understood as an FIO with real phase $\inp{x'-x}{\xi} +
\Delta b_1(z,x',\xi)$ and amplitude $\sigma_A(z,x',\xi) \exp[-\Delta
c_1(z,x',\xi)]$ in $S^m_\rho(X\times \R)$. Applying Theorem
2.5 and the following remark in \cite{Kumano-go:76} we obtain
\addtocounter{delta}{1}
\begin{proposition}
  \label{prop: Hs cont any symbol}
  Let $\A_{(z',z)}$ be the global FIO with kernel
  \begin{align*}
    A_{(z',z)}(x',x) = \int \exp[i
    \phi_{(z',z)}(x',x,\xi)]\  \sigma_A(z,x',\xi)\,
    \dslash \xi
  \end{align*}
  with $\sigma_A(z,.)$ bounded \wrt $z$ with values in
  $S^m(X\times \R^{n})$, $m\in \R$.  Then for all $s
  \in \R$ there exists $M=M(s,m)\geq 0$, $\Delta_\thedelta >0$ such that 
  \begin{align*}
    \|\A_{(z',z)}\|_{(H^{(s)},H^{(s-m)})}
    \leq M\  p(\sigma_A(z,.))
  \end{align*}
  for all $z \in
  [0,Z]$, and $0 \leq\Delta\leq \Delta_\thedelta$, where $p(.)$ is
  some appropriately chosen semi-norm in $S^m(X\times \R^{n})$.
\end{proposition}
This proposition could also be proved by adapting the proof of
  Theorem~\ref{theorem:H^s estimate} to this case. Note that in the
  case $\sigma_A = g_{(z',z)}$ we were able, in the proof of
  Theorem~\ref{theorem:H^s estimate}, to achieve a finer estimate. The
  proof heavily relies on the particular structure of the phase
  function and the amplitude that can be taken as ``close'' as we want
  to those of the identity operator by taking $\Delta$ small enough.

\section{The approximation Ansatz. Convergence in Sobolev spaces}
\label{sec:3}
We first define the Ansatz that approximates the solution operator to
(\ref{eq:one-way})-(\ref{eq:init_cond}). We chose to use a
constant-step subdivision of the interval $[0,Z]$ but the method and
results presented here can be naturally adapted to any subdivision of
$[0,Z]$.
\begin{definition}
  \label{definition:W_P}
  Let $\P=\{z^{(0)},z^{(1)},\dots,z^{(N)}\}$ be a subdivision of
  $[0,Z]$ with $0=z^{(0)}< z^{(1)}<\dots <z^{(N)}=Z$ such that
  $z^{(i+1)} -z^{(i)}=\Delta_\P$.  The operator $\W_{\P,z}$ is defined
  as
  \begin{align*}
    \W_{\P,z} :=
    \left\{
      \begin{array}{ll}
        \G_{(z,0)} & \text{if }\ 0\leq z\leq z^{(1)},\\
        \G_{(z,z^{(k)})} {\displaystyle \prod_{i=1}^{k}} 
	\G_{(z^{(i)},z^{(i-1)})} 
        & \text{if }\ z^{(k)}\leq z\leq z^{(k+1)}.
      \end{array}
    \right.
  \end{align*}
\end{definition}
The following proposition will be useful in the sequel.
\begin{proposition}
  \label{prop:Hs norm under control}
  Let $s \in \R$. There exists $K>0$ such that for every subdivision
  $\P=\{z^{(0)},z^{(1)},\dots,z^{(N)}\}$ of $[0,Z]$ with $0=z^{(0)}<
  z^{(1)}<\dots <z^{(N)}=Z$ and $\W_{\P,z}$ as defined in Definition
  \ref{definition:W_P} we have
  \begin{align*}
    \forall z \in [0,Z], \ \  \| \W_{\P,z} \|_{(H^{(s)},H^{(s)})}\leq K,
  \end{align*}
  if $\Delta_\P$ is small enough. 
\end{proposition}
\begin{proof}
  By Theorem \ref{theorem:H^s estimate} there exits $M>0$ such that if
  $\Delta=z'-z$ is small enough then
  $\|\G_{(z',z)}\|_{(H^{(s)},H^{(s)})} \leq 1+ \Delta M$ for all $z
  \in [0,Z]$; we then obtain $\| \W_{\P,z}\|_{(H^{(s)},H^{(s)})} \leq
  (1+\Delta_\P M)^N = (1+\frac{ZM}{N})^N$ which is bounded as it
  converges to $\exp[ZM]$ as $N$ goes to $\infty$.
\end{proof}

It should be first noticed that $\W_{\P,z}$ is not the solution to
problem (\ref{eq:one-way})-(\ref{eq:init_cond}) even in the case where
the symbols $b$ and $c$ depend only on the transversal variable,
($x$). While singularities propagates along the bicharacteristics
associated with $i a_1 = b_1$, observe however that, with the form of
the phase function $\phi_{(z',z)}$ in (\ref{eq: phase function}), the
operator $\G_{(z',z)}$ propagates singularities along straight lines.
See Part II, for further details, in particular regarding the set
$J_{(z',z)\R}$ that replaces the canonical relation for the
propagation of singularities for FIOs with complex phase
\cite[Sections 25.4-5]{Hoermander:V4}.

Furthermore, by composing the operators $\G_{(z'',z')}$ and
$\G_{(z',z)}$, one convinces oneself that
\begin{align*}
\G_{(z'',z)} \neq \G_{(z'',z')} \circ \G_{(z',z)} 
\end{align*} 
in general if $z''\geq z'\geq z \in [0,Z]$ (use again that
singularities propagate along straight lines). The family of operators
$(\G_{(z',z)})_{(z',z)\in [0,Z]}$ is thus neither a semigroup nor an
evolution system.

We now proceed towards the proof of the convergence of $\W_{\P,z}$ to
the solution operator to problem
(\ref{eq:one-way})-(\ref{eq:init_cond}) in the sense of Sobolev norms
as $N=|\P|$ goes to $\infty$.
\begin{lemma}
  \label{lemma: z' -> Gz'z lipschitz}
  Let $s \in \R$ and $z'',z \in [0,Z]$, with $z <  z''$.  The map $z'
  \mapsto \G_{(z',z)}$, for $z' \in [z,z'']$, is Lipschitz continuous
  with values in $L(H^{(s+1)}(X),H^{(s)}(X))$, for $z''-z=\Delta$
  small enough. More precisely there exists $C>0$ such that for all
  $u_0 \in H^{(s+1)}(X)$ and $z^{(1)}, z^{(2)} \in [z,z'']$
  \begin{align}
    \label{eq: z' -> Gz'z Lipschitz}
    \|(\G_{(z^{(2)},z)}-\G_{(z^{(1)},z)})(u_0)\|_{H^{(s)}} 
    \leq C |z^{(2)}-z^{(1)}| \|u_0\|_{H^{(s+1)}}.
  \end{align}
\end{lemma}
\begin{proof}
  Let $z^{(1)},z^{(2)} \in [z'',z]$ and let $u_0 \in H^{s+1}(X)$.
  Write 
  \begin{multline*}
    (\G_{(z^{(2)},z)}-\G_{(z^{(1)},z)})(u_0) (x') =\\
    - \int_{z^{(1)}}^{z^{(2)}}\iint \exp[i \inp{x'-x}{\xi}
    - (z'-z)a(z,x',\xi)]\
        a(z,x',\xi)\ \ u_0(x)\ d x\ \dslash \xi \ d z'.
  \end{multline*}
  When $\Delta$ is small enough we can apply Proposition~\ref{prop: Hs
    cont any symbol} and  obtain (\ref{eq: z' -> Gz'z Lipschitz})
\end{proof}
\begin{lemma}
  \label{lemma: cont G, dz G}
  Let $s \in \R$, $z'',z \in [0,Z]$, with $z <  z''$, and let $u_0
  \in H^{(s+1)}(X)$.  Then the map $z'\mapsto \G_{(z',z)}(u_0)$ is in
  $\Con^0([z,z''],H^{(s+1)}(X)) \cap \Con^1([z,z''],H^{(s)}(X))$ for
  $z''-z=\Delta$ small enough.
\end{lemma}
\begin{proof}
  Let $z^{(1)} \in [z,z'']$ and let $\varepsilon>0$.  Choose $z''-z$
  small enough such that Theorem~\ref{theorem:H^s estimate} and
  Lemma~\ref{lemma: z' -> Gz'z lipschitz} apply and Choose $u_1 \in
  H^{(s+2)}$ such that $\| u_0 - u_1\|_{H^{(s+1)}}\leq \varepsilon$.
  Then for $z^{(2)} \in [z,z'']$
  \begin{multline}
    \label{estimate: Hs+1 cont}
    \| \G_{(z^{(2)},z)}(u_0) - \G_{(z^{(1)},z)}(u_0)\|_{H^{(s+1)}} \leq 
    \|\G_{(z^{(2)},z)}(u_0-u_1)\|_{H^{(s+1)}} \\
    + \|\G_{(z^{(2)},z)}(u_1) - \G_{(z^{(1)},z)}(u_1)\|_{H^{(s+1)}}
    + \|\G_{(z^{(1)},z)}(u_0-u_1)\|_{H^{(s+1)}}\\
    \leq 2 (1+ \Delta M) \varepsilon + C|z^{(2)}-z^{(1)}| \| u_1\|_{H^{(s+2)}}.
    \end{multline}
    The continuity of the map follows. Differentiating
    $G_{(z',z)}(u_0)$ \wrt $z'$ we can prove that the resulting map
    $z' \mapsto \d_{z'}G_{(z',z)}(u_0)$ is Lipschitz continuous with
    values in $L(H^{(s+2)},H^{(s)})$ following the proof of
    Lemma~\ref{lemma: z' -> Gz'z lipschitz}: there exists $C>0$ such
    that for all $v \in H^{(s+2)}(X)$
    \begin{align*}
      \label{eq: z' -> d_z' Gz'z Lipschitz}
      \|(\d_{z'} \G_{(z^{(2)},z)}-\d_{z'}
      \G_{(z^{(1)},z)})(v)\|_{H^{(s)}} \leq C |z^{(2)}-z^{(1)}|
      \|v\|_{H^{(s+2)}}.
    \end{align*} 
    We also see that the map $v \mapsto \d_{z'}G_{(z',z)}(v)$ is
    continuous from $H^{(s+1)}$ into $H^{(s)}$ with bounded continuity
    module according to Proposition~\ref{prop: Hs cont any symbol}.
    With $u_0 \in H^{(s+1)}(X)$ we make a similar choice for $u_1
    \in H^{(s+2)}(X)$ and  obtain an estimate for
    \[ \| \d_{z'} \G_{(z^{(2)},z)}(u_0) - \d_{z'}
    \G_{(z^{(1)},z)}(u_0)\|_{H^{(s)}}
    \]
    of the same form as in (\ref{estimate: Hs+1 cont}).
 \end{proof}
The two previous lemmas yield
\begin{proposition} 
  Let $s \in \R$, $\P$ a subdivision of $[0,Z]$ as in
  Definition~\ref{definition:W_P} and let $u_0 \in H^{(s+1)}(X)$.
  Then the map $\W_{\P,z} (u_0)$ is $\Con^0([0,Z],H^{(s+1)}(X))$
  and piecewise $\Con^1([0,Z],H^{(s)}(X))$ if $\P$ is chosen such
  that $\Delta_\P$ is small enough.  The map $z\mapsto \W_{\P,z}
  (u_0)$ is in fact globally Lipschitz with $C>0$ such that
  \[
   \| \W_{\P,z'} (u_0) - \W_{\P,z} (u_0)\|_{H^{(s)}} 
   \leq C |z'-z| \| u_0 \|_{H^{(s+1)}}.
  \]
\end{proposition}
We recall that $U(z',z)$ is the solution operator of the Cauchy
problem (\ref{eq:one-way})-(\ref{eq:init_cond}). We can then apply the
energy estimate~(\ref{eq:energy estimate}) to $U(z,0)(u_0) -
\W_{\P,z}(u_0)$ (adapt the proof of Lemma 23.1.1 in
\cite{Hoermander:V3} to the case of a Lipschitz piecewise $C^1$
function) and obtain
\begin{multline}
  \label{eq: energy estimate2}
  \sup_{z\in [0,Z]} \exp [-\lambda z]\ \|U(z,0)(u_0) -
  \W_{\P,z}(u_0)\|_{H^{(s)}}
  \\
  \leq 2 \int_0^Z \exp[-\lambda z]\ \|(\d_z + a_z(x,D_x))
  \W_{\P,z}(u_0)\|_{H^{(s)}} d z.
\end{multline}
Let $u_0 \in H^{(s+1)}(X)$ and let $\P= \{z^{(0)}, \dots, z^{(N)} \}$.  We
take $z \in ]z^{(k)},z^{(k+1)}[$.  Then
\begin{align*}
  (\d_z + a_z(x,D_x))\ & \W_{\P,z}(u_0)\\
 &= (\d_z + a_z(x,D_x)) \left( \G_{(z,z^{(k)})} 
 {\displaystyle \prod_{i=1}^{k}} \G_{(z^{(i)},z^{(i-1)})} (u_0)\right) \\
 &= (\d_z + a_z(x,D_x)) \left(\G_{(z,z^{(k)})} (u_k) \right)
\end{align*}
with $u_k := {\displaystyle \prod_{i=1}^{k}} \G_{(z^{(i)},z^{(i-1)})} (u_0)$
which is in $H^{(s+1)}(X)$ by Theorem~\ref{theorem:H^s
  estimate}.  We first turn our attention towards the term $(\d_z +
a_z(x,D_x)) \left( \G_{(z,z^{(k)})} (u) \right)$ for any $u \in
H^{(s+1)}(X)$ as the norm of $u_k$ in $H^{(s+1)}(X)$ remains
under control even if $|\P|=N$ becomes very large by
Proposition~\ref{prop:Hs norm under control}:
\begin{multline}
  \label{eq: HS norm under control}
  \exists K>0,\ \ \| u_k \|_{H^{(s+1)}} \leq K \| u_0 \|_{H^{(s+1)}}, \  
  k \in \{0,\dots,N \},\\
  N=|\P|  \in \N, \  u_0 \in H^{(s+1)}(X),
\end{multline}
if $\Delta_\P$ is small enough.

We shall need the following lemma which is a variant to Lemma~\ref{lemma:symb - amp}
\begin{lemma}
  \label{lemma:symb - amp 2}
  Let $q_\Delta(z,x,y,\xi)$ be an amplitude in
  $S^0_\rho(\R^{2p}\times \R^p)$ depending on the parameters
  $\Delta\geq 0$ and $z \in [0,Z]$ that satisfies Property
  (\ref{eq:Q_L}) and such that $q_\Delta(z,.)|_{\Delta=0}=0$. Let
  $r(x,y,\xi) \in S^s(\R^{2p}\times \R^p)$ for some $s \in
  \R$.  Then
  \begin{align*}
    \symb{q_\Delta\ r } (z,x,\xi) - q_\Delta(z,x,x,\xi)\ r(x,x,\xi) =
    \Delta^{m + 2 \delta} \lambda_{\Delta }^{m}(z,x,\xi), \ 0\leq m \leq \rho
    - \delta,
  \end{align*}
  where the function $\lambda_\Delta^{m }(z,x,\xi)$ is bounded with
  respect to $\Delta$ and $z$ with values in $S^{m +s - (\rho -
    \delta)}_\rho(\R^p\times \R^p)$.
\end{lemma}
\begin{proof}
  We proceed as in the proof of Lemma~\ref{lemma:symb - amp} (we take
  $p=1$ for the sake of concision). We obtain
  \begin{multline*}
    \symb{q_\Delta r} (z,x,\xi) - q_\Delta(z,x,x,\xi)r(x,x,\xi) \\
    = \symb{i \int_0^1 \d_3 \d_4 \lambda_\Delta(z,x,(1-s)x+s y,\xi)\ ds},
  \end{multline*}
  where here 
  \begin{multline*}
    \d_3 \d_4 \lambda_\Delta(z,x,y,\xi) = (\d_y \d_\xi
    q_\Delta)(z,x,y,\xi) \ r(x,y,\xi) + \d_y q_\Delta(z,x,y,\xi)
    \d_\xi r(x,y,\xi) \\
    + \d_\xi q_\Delta(z,x,y,\xi) \d_y r(x,y,\xi) +
    q_\Delta(z,x,y,\xi) \d_y\d_\xi r(x,y,\xi).
  \end{multline*}
  The first two terms are treated like in the proof of
  Lemma~\ref{lemma:symb - amp}. For the Third term, with Property
  (\ref{eq:Q_L}) we write
  \begin{align*}
    \d_\xi q_\Delta\ \d_y r = \Delta^{m'+\delta}\ q_\Delta^{m'(00)1}\ 
    \d_y r, \ \ 0\leq m' \leq 1-\delta
  \end{align*}
  where $q_\Delta^{m'(0,0)1}\ d_y r \in
  S^{m'+s-\rho}_\rho(\R^{2p}\times \R^p)$. We actually take
  $\delta \leq m' \leq 1-\delta$ and write $m=m'-\delta$.  We obtain
  \begin{align*}
    \d_\xi q_\Delta\  \d_y r 
    = \Delta^{m+2\delta}\  \tilde{q}_\Delta^{m},
  \end{align*}
  where $\tilde{q}_\Delta^{m}$ is bounded \wrt $\Delta$ with values
  in $S^{m+s-\rho+\delta}_\rho(\R^{2p}\times \R^p)$ and
  $0\leq m\leq 1-2\delta= \rho -\delta$.
  For the fourth term we write 
  \begin{align*}
    q_\Delta = \Delta^{m'} q_\Delta^{m'(0,0)0}, \ \ 0\leq m' \leq 1, 
  \end{align*}
  where $q_\Delta^{m'(0,0)0} \in S^{m'}_\rho(\R^{2p}\times
  \R^p) $ by Property (\ref{eq:Q_L}) since $q_\Delta|_{\Delta=0}=0$.
  We actually take $2\delta \leq m' \leq 1$ and write $m=m'-2\delta$.
  Then
  \begin{align*}
    q_\Delta \d_y\d_\xi r = \Delta^{m+2\delta} \hat{q}_\Delta^{m},
  \end{align*}
  where $\hat{q}_\Delta^{m}$ is bounded \wrt $\Delta$ with values in
  $S^{m+s-(\rho-\delta)}_\rho(\R^{2p}\times \R^p)$ as
  $m+s-1+2\delta =m+s-(\rho-\delta)$ and $0\leq m\leq 1-2\delta= \rho
  -\delta$. We conclude like in the proof of Lemma~\ref{lemma:symb -
    amp}.
\end{proof}

For the next proposition we shall need the following assumption as
announced in Section~\ref{sec:1}
\begin{asmpt}
  \label{assumpt:Lipschitz assumption2}
  The symbol $a(z,.)$ is assumed to be in ${\mathscr L}([0,Z],S^1(\R^{n}\times
  \R^{n}))$, i.e. Lipschitz continuous \wrt $z$ with values in
  $S^1(\R^{n}\times \R^{n})$, in the sense that,
\begin{align*}
  a(z',x,\xi) - a(z,x,\xi) = (z'-z) \tilde{a}(z',z,x,\xi),\  \ 
  0 \leq z \leq z' \leq Z
\end{align*}
  with $\tilde{a}(z',z,x,\xi)$ bounded \wrt $z'$ and $z$ with values
in $S^1(\R^{n}\times \R^{n})$.
\end{asmpt}
\addtocounter{delta}{1}
\begin{proposition}
  Let $s \in \R$.  There exists $\Delta_\thedelta > 0$
  and $C\geq 0$ such that for $z'-z=\Delta$, $\Delta \in
  [0,\Delta_\thedelta]$,
\[
\| (\d_{z'} + a_{z'}(x,D_x))
 \G_{(z',z)} \|_{(H^{(s)}, H^{(s-1)})} 
  \leq C \Delta^{\hf} 
\]
\end{proposition}
Like in the proof of Theorem~\ref{theorem:H^s estimate} we assume that
$c_1$ satisfies property $P_L$ for some $L \geq 2$. We know that it is
always true for $L=2$ by Lemma~\ref{lemma: PL with L=2 for c_1} but
special choices for $c_1$ can be made. As before we use $\rho=1-1/L$
and $\delta=1/L$ with $\rho > \delta$ for $L>2$ and $\rho=\delta=\hf$
for $L=2$.
\begin{proof}
  Let $\A_{(z',z)}$ be $\d_{z'} \G_{(z',z)}$ and $\B_{(z',z)}$ be
  $a_{z'}(x,D_x) \circ \G_{(z',z)}$ with respective kernels
  $A_{(z',z)}(x',x)$ and $B_{(z',z)}(x',x)$. We have
  \begin{align*}
    A_{(z',z)} (x',x) = - \int \exp[i \inp{x'-x}{\xi}]
    \exp[ - \Delta a(z,x',\xi)]\: 
     a(z,x',\xi) \ \dslash \xi.
   \end{align*}
   Let us define
  \[
  {\cal D}_{(z',z)}: = (\A_{(z',z)} + \B_{(z',z)}) \circ
  E^{-2s}\circ (\A_{(z',z)} + \B_{(z',z)})^\ast. 
  \]
  We prove in the following lemma that for $r,s \in \R$, $\| {\cal
    D}_{(z',z)} \|_{(H^{(r)},H^{(r+2s-2))})} \leq C \Delta$
  uniformly \wrt $z \in [0,Z]$ for $\Delta$ small enough.
  The conclusion then follows: if ${\cal C}_{(z',z)}:=
  E^{s-1}\circ {\cal D}_{(z',z)} \circ E^{s-1}$ then
  $\| {\cal C}_{(z',z)} \|_{(L^2,L^2)}\leq C \Delta$ (take $r =
  -s+1$); then $\| E^{s-1} \circ (\A_{(z',z)} +
  \B_{(z',z)})\circ E^{-s}\|_{(L^2,L^2)}\leq C
  \Delta^{\hf}$.
\end{proof}
\begin{lemma}
  \label{lemma: Delta coming out of symbol}
  Let $r,s \in \R$. Then $\| {\cal D}_{(z',z)}
  \|_{(H^{(r)},H^{(r+2s-2))})} \leq C \Delta$ uniformly \wrt $z
  \in [0,Z]$ for $\Delta$ small enough.
\end{lemma}
\begin{proof}
  The operator ${\cal D}_{(z',z)}$ is made up of four terms:
  \begin{align*}
    & {\cal D}_{1,(z',z)}:= \A_{(z',z)} \circ
    E^{-2s}\circ \A_{(z',z)}^\ast, \ \ \ 
    {\cal D}_{2,(z',z)}:= \A_{(z',z)} \circ
    E^{-2s}\circ \B_{(z',z)}^\ast,\\
    & {\cal D}_{3,(z',z)}:= \B_{(z',z)} \circ
    E^{-2s}\circ \A_{(z',z)}^\ast, \ \ \ 
    {\cal D}_{4,(z',z)}:= \B_{(z',z)} \circ
    E^{-2s}\circ \B_{(z',z)}^\ast.
  \end{align*}
The kernel of ${\cal D}_{1,(z',z)}$ is given by 
\begin{multline*}
  D_{1,(z',z)}(x',x) \\
  = \int \exp\left[ i\inp{x'-x}{\xi}   +
      i\Delta \left(b_1(z,x',\xi)-b_1(z,x,\xi)\right)\right]\
    \tilde{d}_{1,z}(x',x,\xi)\
    \dslash \xi
\end{multline*}
where
\begin{align*}
  \tilde{d}_{1,z}(x',x,\xi) 
  = \omega_{(z',z)}(x',x,\xi) \ 
  a(z,x',\xi) \ 
  \overline{a}(z,x,\xi), 
\end{align*}
and
\begin{multline*}
 \omega_{(z',z)}(x',x,\xi):= g_{(z',z)}(x',\xi) \ 
  \overline{g_{(z',z)}}(x,\xi)
\exp[-\Delta (c_1(z,x',\xi)+c_1(z,x,\xi))]\\
\e{\xi}{-2s}
\end{multline*} 
with $g_{(z',z)}$ given in (\ref{eq: amplitude g}).
Following the proof of Theorem~\ref{theorem:H^s estimate} we write
$b_1(z,x',\xi)-b_1(z,x,\xi) = \inp{x'-x}{ h(z,x',x,\xi)}$ where $h$ is
homogeneous of degree one in $\xi$,
  $|\xi|\geq 1$.  The function $h$ and continuous \wrt $z$
  with values in $S^1(X\times \R^n)$. We thus obtain that the change of
variables $\xi \to \xi+\Delta h(z,x',x,\xi)$ is a global diffeomorphism
for $\Delta$ small enough (uniformly in $z \in [0,Z]$).  The Jacobian
$\J_\Delta(z,x',x,\xi)$ is homogeneous of degree zero in $\xi$,
$\Cinf$ \wrt $\Delta$ and bounded \wrt $z$
with values in $S^0(\R^{2n}\times \R^{n})$.  We then
have
\begin{multline*}
  D_{1,(z',z)}(x',x) =\! \int \exp\left[ i\inp{x'-x}{\xi} \right]
  \tilde{d}_{1,z}(x',x,\tilde{\xi}(\Delta,\xi))\ \J(\Delta,z,x',x,\xi)\  
  \dslash \xi.
\end{multline*}
The function $\tilde{\xi}(\Delta,z,x',x,\xi)$, written
$\tilde{\xi}(\Delta,\xi)$ for concision, is bounded \wrt  $z$ and
$\Cinf$ \wrt $\Delta$ in $S^1(\R^{2n}\times\R^{n})$ and homogeneous
of degree 1 in $\xi$ as shown in Lemma~\ref{lemma: change of
  variables}.  It follows that
$\tilde{d}_{1,z}(x',x,\tilde{\xi}(\Delta,\xi))\ \J(\Delta,z,x',x,\xi)$
is then bounded \wrt $z$ and $\Delta$ with values in
$S_\rho^{2-2s}(\R^{2n} \times \R^n)$ by Lemma~\ref{lemma:
  change of variable - symbol} and the proof of
Theorem~\ref{theorem:H^s estimate}.  Note that if $\Delta=0$ then
$\tilde{\xi}(\Delta,\xi) = \xi$. The operator ${\cal D}_{1,(z',z)}$ is
thus in $\Psi^{2-2s}_\rho$ with symbol
\[
d_{1,(z',z)}(x',\xi)=
\symb{\tilde{d}_{1,z}(x',x,\tilde{\xi}(\Delta,\xi))\ 
  \J(\Delta,z,x',x,\xi)}(x',\xi).
\] 
Similarly we prove that $\A_{(z',z)} \circ E^{-2s}\circ
\G_{(z',z)}^\ast$ is the $\psi$DO with {\em amplitude}
\[
  - \omega_{(z',z)}(x',x,\tilde{\xi}(\Delta,\xi)) \ 
  a(z,x',\tilde{\xi}(\Delta,\xi))\  \J(\Delta,z,x',x,\xi).
\]
The operator ${\cal D}_{2,(z',z)}$ is thus in $\Psi^{2-2s}_\rho(X)$ with  symbol
\begin{multline*}
  d_{2,(z',z)}(x',\xi) \\= -
  \symb{\omega_{(z',z)}(x',x,\tilde{\xi}(\Delta,\xi))\ 
  a(z,x',\tilde{\xi}(\Delta,\xi)) 
  \J(\Delta,z,x',x,\xi)}\ \#\ a^\ast(z',x',\xi)
\end{multline*}
Similarly we find that the operators ${\cal D}_{3,(z',z)}$ and ${\cal
  D}_{4,(z',z)}$ are in $\Psi^{2-2s}_\rho(X)$ with respective
symbols
\begin{multline*}
  d_{3,(z',z)}(x',\xi) = - a(z',x',\xi)\\
  \#\ \symb{\omega_{(z',z)}(x',x,\tilde{\xi}(\Delta,\xi))\ 
    \J(\Delta,z,x',x,\xi)\ \overline{a}(z,x,\tilde{\xi}(\Delta,\xi))}
\end{multline*}
and 
\begin{multline*}
  d_{4,(z',z)}(x',\xi) = a(z',x',\xi)\ \#\ 
  \symb{\omega_{(z',z)}(x',x,\tilde{\xi}(\Delta,\xi))\ 
    \J(\Delta,z,x',x,\xi)}\\ \#\ a^\ast(z',x',\xi),
\end{multline*}
For $q(x',x,\xi)$ an amplitude we define
\begin{multline*}
  \Sigma\{q\}(x',\xi) :=\sigma\{ \e{\xi}{-2s} \ a(z,x',\xi)\ 
  q(x',x,\xi)\ 
  \overline{a}(z,x,\xi)\}\\
  - \sigma\{\e{\xi}{-2s} \ a(z,x',\xi)\ q(x',x,\xi) \}\ \#\  a^\ast(z,x',\xi)\\
  + a(z,x',\xi)\ \#\ \sigma\{\e{\xi}{-2s}\ q(x'x,\xi)\}
  \#\  a^\ast(z,x',\xi)\\
  - a(z,x',\xi)\ \#\ \sigma\{ \e{\xi}{-2s}\ 
  q(x'x,\xi) \ \overline{a}(z,x,\xi)\}.
\end{multline*}
The operator ${\cal D}_{(z',z)}$ is thus in $\Psi^{2-2s}_\rho(X)$
with symbol
\begin{align*}
  d_{(z',z)} = d_{1,(z',z)}+ d_{2,(z',z)}+ d_{3,(z',z)}+ d_{4,(z',z)}.
\end{align*}
Such a symbol is bounded \wrt $\Delta$, for $\Delta$ small
enough, as the composition formula for symbols is a bounded map.
Note that 
\begin{align*}
  g_{(z',z)}(x',\xi) \ \overline{g_{(z',z)}}(x,\xi) \e{\xi}{-2s}
  \J(\Delta,z,x',x,\xi) = \e{\xi}{-2s} + \Delta k_\Delta (z,x',x,\xi)
\end{align*}
with $k_\Delta$ bounded \wrt $z$ and $\Cinf$ \wrt $\Delta$ with values
in $S^{-2s}(X'\times X \times \R^n)$ as $\omega_{(z',z)}\J$ is itself
$\Cinf$ \wrt $\Delta$ and $g_{(z',z)}\ \overline{g_{(z',z)}}
\J|_{\Delta=0}=1$.  By Assumption~\ref{assumpt:Lipschitz assumption2},
we also write $a(z',x,\xi) = a(z,x,\xi) + \Delta
\tilde{a}(z',z,x,\xi)$ with $\tilde{a}(z',z,.)$ bounded \wrt $z$ and
$\Delta$ with values in $S^1(X\times \R)$.  We thus obtain
\begin{align*}
 {\cal D}_{(z',z)} &= {\cal D}^a_{(z',z)} + \Delta {\cal D}^1_{(z',z)}
\end{align*}
with symbols
\begin{align*}
 d^a_{(z',z)}&:= \Sigma\{ \e{\xi}{-2s} p_\Delta (z,x',x,\xi) \}, \\ 
\end{align*}
and $d^1_{(z',z)}$ which is bounded \wrt $z$ and $\Delta$ with values
in $S^{2-2s}_\rho (X'\times X\times \R^n)$. The
Calder\'{o}n-Vaillancourt theorem (see \cite[Chapter 7, Sections
1,2]{Kumano-go:81} or \cite[Section XIII-2]{Taylor:81}) in the case
$L=2$ or Theorem 18.1.11 in \cite{Hoermander:V3} in the case $L>2$
yields $\|{\cal D}^1_{(z',z)}\|_{(H^{(r)},H^{(r+2s-2))})} \leq K_1$.
Note that for a symbol $q(x',\xi)$ we have $\Sigma\{q(x',\xi)\} =0$ as
\begin{align*}
  \sigma\{ q(x',\xi) \ 
  \overline{a}(z,x,\xi)\} =  q(x',\xi)\ \#\  a^\ast(z,x',\xi)
  = \sigma\{q(x',\xi) \}\ \#\ a^\ast(z,x',\xi),
\end{align*}
for any symbol $q$.
Thus $d^a{(z',z)}= \Sigma\{ \e{\xi}{-2s} (p_\Delta (z,x',x,\xi)-1) \}$.
Lemma~\ref{lemma:symb - amp 2} allows us to  write (take $m=\rho- \delta$) 
\begin{multline*}
  \sigma\{\e{\xi}{-2s} (p_\Delta (z,x',x,\xi)-1) 
  a(z,x',\xi)\overline{a}(z,x,\xi)\}\\
  = \e{\xi}{-2s} (p_\Delta (z,x',x',\xi)-1)a(z,x',\xi)  \overline{a}(z,x',\xi)
  + \Delta \lambda_{\Delta,1} (z,x',\xi)
\end{multline*}
where $\lambda_{\Delta,1}$ bounded \wrt $z$ and $\Delta$ with values
in $S^{2-2s}_\rho(X' \times\R^n)$.
We also write 
\begin{multline*}
  \sigma\{\e{\xi}{-2s} (p_\Delta (z,x',x,\xi)-1) 
  a(z,x',\xi)\}\ \# \ a^\ast (z,x',\xi)\\
  = (\e{\xi}{-2s} (p_\Delta (z,x',x',\xi)-1)a(z,x',\xi))\ 
  \# \ a^\ast (z,x',\xi)\\
  + \Delta \lambda_{\Delta,2} (z,x',\xi)\ \# \ a^\ast (z,x',\xi)\\
  = \sigma\{\e{\xi}{-2s} (p_\Delta (z,x',x',\xi)-1)a(z,x',\xi) 
  \overline{a} (z,x,\xi)\} \\
  + \Delta \lambda_{\Delta,2} (z,x',\xi)\ \# \ a^\ast (z,x',\xi)\\
  = \e{\xi}{-2s} (p_\Delta (z,x',x',\xi)-1)a(z,x',\xi) 
  \overline{a} (z,x',\xi)\\
  + \Delta (\lambda_{\Delta,3} (z,x',\xi) +
  \lambda_{\Delta,2} (z,x',\xi)\ \# \ a^\ast (z,x',\xi))
\end{multline*}
where $\lambda_{\Delta,2}$ and $\lambda_{\Delta,3}$ are bounded \wrt
$z$ and $\Delta$ with values in $S^{1-2s}_\rho(X' \times\R^n)$
and $S^{2-2s}_\rho(X' \times\R^n)$ respectively.  Similarly we
have
\begin{multline*}
  \sigma\{\e{\xi}{-2s} (p_\Delta (z,x',x,\xi)-1) \}\ \# \ a^\ast (z,x',\xi)\\
  = (\e{\xi}{-2s} (p_\Delta (z,x',x',\xi)-1))\ 
  \# \ a^\ast (z,x',\xi)
  + \Delta \lambda_{\Delta,4} (z,x',\xi)\ \# \ a^\ast (z,x',\xi)\\
  = \sigma\{\e{\xi}{-2s} (p_\Delta (z,x',x',\xi)-1) 
  \overline{a} (z,x,\xi)\} 
  + \Delta \lambda_{\Delta,4} (z,x',\xi)\ \# \ a^\ast (z,x',\xi)\\
  = \e{\xi}{-2s} (p_\Delta (z,x',x',\xi)-1)
  \overline{a} (z,x',\xi)\\
  + \Delta (\lambda_{\Delta,5} (z,x',\xi) +
  \lambda_{\Delta,4} (z,x',\xi)\ \# \ a^\ast (z,x',\xi))
\end{multline*}
where $\lambda_{\Delta,4}$ and $\lambda_{\Delta,5}$ are bounded \wrt
$z$ and $\Delta$ with values in $S^{-2s}_\rho(X'
\times\R^n)$ and $S^{1-2s}_\rho(X' \times\R^n)$
respectively and 
\begin{multline*}
  \sigma\{\e{\xi}{-2s} (p_\Delta (z,x',x,\xi)-1) 
  \overline{a}(z,x,\xi)\}\\
  = \e{\xi}{-2s} (p_\Delta (z,x',x',\xi)-1) \overline{a}(z,x',\xi)
  + \Delta \lambda_{\Delta,6} (z,x',\xi)
\end{multline*}
where $\lambda_{\Delta,6}$ bounded \wrt $z$ and $\Delta$ with values
in $S^{1-2s}_\rho(X' \times\R^n)$.
We thus obtain
\begin{align*} 
d^a_{(z',z)} = \Delta (\lambda_{\Delta,1} + \lambda_{\Delta,3} 
+ \lambda_{\Delta,2}\ \# \ a^\ast 
+ a\ \# \ \lambda_{\Delta,5}  
+ a\ \#\ \lambda_{\Delta,4}\ \# \ a^\ast) = \Delta \tilde{d}^a_{(z',z)}
\end{align*}
with $\tilde{d}^a_{(z',z)}$ bounded \wrt  $z$ and $\Delta$ with
values in $S^{2-2s}_\rho(X' \times\R^n)$. This
concludes the proof.
\end{proof}
We have thus obtained a convergence result in the Sobolev space
$H^{(s)}(\R^n)$ for $\W_{\P,z}(u_0)$ if the initial data $u_0$ is in
$H^{(s+1)}(\R^n)$. The result is actually the convergence of the Ansatz
$\W_{\P,z}$ to the solution operator $U(z,0)$ in the norm of
$L(H^{(s+1)}(\R^n),H^{(s)}(\R^n))$:
\begin{theorem}
  \label{theorem: convergence result}
  Assume that $a(z,.)$ is in ${\mathscr L}([0,Z],S^1(\R^{n}\times
  \R^{n}))$, i.e. Lipschitz continuous \wrt $z$ with values in
  $S^1(\R^{n}\times \R^{n})$, in the sense that,
  \begin{align*}
    a(z',x,\xi) - a(z,x,\xi) = (z'-z) \tilde{a}(z',z,x,\xi),\  \ 
    0 \leq z \leq z' \leq Z
  \end{align*}
  with $\tilde{a}(z',z,x,\xi)$ bounded \wrt $z'$ and $z$ with values
  in $S^1(\R^{n}\times \R^{n})$.  Let $s \in \R$. Then the
  approximation Ansatz $\W_{\P,z}$ converges to the solution operator
  $U(z,0)$ of the Cauchy
  problem~(\ref{eq:one-way})-(\ref{eq:init_cond}) in
  $L(H^{(s+1)}(\R^n),H^{(s)}(\R^n))$ as $\Delta_\P$ goes to 0 with a
  convergence rate of order $\hf$:
  \begin{align*}  
    \| \W_{\P,z}  - U(z,0) \|_{(H^{(s+1)},H^{(s)})} \leq C \Delta_\P^\hf.
  \end{align*}
\end{theorem}
\begin{proof}
  Using (\ref{eq: energy estimate2}) and (\ref{eq: HS norm under control}) we obtain
  \begin{multline*}
    \sup_{z\in [0,Z]} \exp [-\lambda z]\ \|U(z,0)(u_0) -
    \W_{\P,z}(u_0)\|_{H^{(s)}}
    \\
    \leq 2 \int_0^Z \exp[-\lambda z]\ \Delta_\P^{\hf} C K 
    \| u_0\|_{H^{(s+1)}} d z \leq C \Delta_\P^{\hf} \| u_0\|_{H^{(s+1)}}.
\end{multline*}
The result follows.
\end{proof}

If we change the assumption made on the symbol $a(z,.)$ to some H\"older
type continuity, then the corresponding change in the proof of
Lemma~\ref {lemma: Delta coming out of symbol} yields the following
weaker result
\begin{theorem}
  \label{theorem: convergence result hoelder}
  Assume that $a(z,.)$ is in $\Con_\alpha([0,Z],S^1(\R^{n}\times
  \R^{n}))$, i.e. H\"older  continuous \wrt $z$ with values in
  $S^1(\R^{n}\times \R^{n})$, in the sense that, for some $0< \alpha <1$
  \begin{align*}
    a(z',x,\xi) - a(z,x,\xi) = (z'-z)^\alpha\  \tilde{a}(z',z,x,\xi),\  \ 
    0 \leq z \leq z' \leq Z
  \end{align*}
  with $\tilde{a}(z',z,x,\xi)$ bounded \wrt $z'$ and $z$ with values
  in $S^1(\R^{n}\times \R^{n})$.
  Let $s \in \R$.  Then the
  approximation Ansatz $\W_{\P,z}$ converges to the solution operator
  $U(z,0)$ of the Cauchy
  problem~(\ref{eq:one-way})-(\ref{eq:init_cond}) in
  $L(H^{(s+1)}(\R^n),H^{(s)}(\R^n))$ as $\Delta_\P$ goes to 0 with a
  convergence rate of order $\alpha/2$:
  \begin{align*}  
    \| \W_{\P,z}  - U(z,0) \|_{(H^{(s+1)},H^{(s)})} 
    \leq C \Delta_\P^{\alpha/2}.
  \end{align*}
\end{theorem}

A result similar to that of the previous theorems can be obtained with
weaker assumptions, namely without assumptions on the symbol $a(z,.)$
like those made in Theorems~\ref{theorem: convergence result} and
\ref{theorem: convergence result hoelder}, by introducing another,
yet natural, Ansatz to approximate the exact solution to the Cauchy
problem (\ref{eq:one-way})-(\ref{eq:init_cond}).  For a symbol
$q(z,y,\eta) \in \Con^0([0,Z],S^m(\R^{p}\times \R^{r}))$ we define
$\widehat{q}_{(z',z)}(y,\eta) \in \Con^0([0,Z]^2,S^m(\R^{p}\times
\R^{r}))$
\begin{align*}
  & \widehat{q}_{(z',z)}(y,\eta) := \frac{1}{z'-z}\int_z^{z'} q(s,y,\eta)\ d s.\ 
\end{align*}
Then we define
\begin{multline}
  \label{eq: phase function new}
  \widehat{\phi}_{(z',z)}(x',x,\xi):= \inp{x'-x}{\xi}
  +i \Delta \widehat{a}_{1(z',z)}(x',\xi)\\
  = \inp{x'-x}{\xi}  + \Delta \widehat{b}_{1(z',z)}(x',\xi) +i
  \Delta \widehat{c}_{1(z',z)}(x',\xi).
\end{multline}
and 
\begin{align}
   \label{eq: amplitude g new}
   \widehat{g}_{(z',z)}(x,\xi) := \exp[-\Delta \widehat{a}_{0(z',z)}(x,\xi)].
\end{align}
and finally, following \cite{N.Kumano-go:95}, we denote
$\widehat{\G}_{(z',z)}$ the FIO with distribution kernel
\begin{multline*}
  \widehat{G}_{(z',z)}(x',x) = \int \exp[i \inp{x'-x}{\xi} ]
  \exp[ - \Delta \widehat{a}_{(z',z)}(x',\xi)]\, \dslash \xi\\
  = \int \exp[i \widehat{\phi}_{(z',z)}(x',x,\xi)]\: 
  \widehat{g}_{(z',z)}(x',\xi) \, \dslash \xi.
\end{multline*}
with the associated approximation Ansatz
\begin{definition}
  \label{definition:W_P new}
  Let $\P=\{z^{(0)},z^{(1)},\dots,z^{(N)}\}$ be a subdivision of
  $[0,Z]$ with $0=z^{(0)}< z^{(1)}<\dots <z^{(N)}=Z$ such that
  $z^{(i+1)} -z^{(i)}=\Delta_\P$.  The operator $\widehat{\W}_{\P,z}$ is defined
  as
  \begin{align*}
    \widehat{\W}_{\P,z} :=
    \left\{
      \begin{array}{ll}
        \widehat{\G}_{(z,0)} & \text{if }\ 0\leq z\leq z^{(1)},\\
        \widehat{\G}_{(z,z^{(k)})} {\displaystyle \prod_{i=1}^{k}} 
	\widehat{\G}_{(z^{(i)},z^{(i-1)})} 
        & \text{if }\ z^{(k)}\leq z\leq z^{(k+1)}.
      \end{array}
    \right.
  \end{align*}
\end{definition}
Most results of Sections~\ref{sec:2} and ~\ref{sec:3} apply to this
new Ansatz. We give some details about how to adapt some of the proofs.
We have
\begin{lemma}
  Let $q(z,y,\eta) \in \Con^0([0,Z],S^1(\R^{p}\times \R^{r}))$ that
  satisfies Property (\ref{eq:P_L}). Then $\widehat{q}_{(z',z)}(y,\eta)$
  also satisfies Property (\ref{eq:P_L}).
\end{lemma}
Property (\ref{eq:P_L}) in Definition~\ref{definition:special damping}
is now to be understood \wrt to two parameters $z'$ and $z$.
\begin{proof}
  Uniform bounds \wrt $z$ and $z'$ will be immediate.
  The case $|\alpha| + |\beta|\geq L$ is clear by Remark~\ref{rem:
  P_L for alpha + beta >= L}. Let then $|\alpha| + |\beta| < L$ and
  observe that 
  \begin{multline*}
    |\d_y^\alpha d_\eta^\beta \widehat{q}_{(z',z)}(y,\eta) |
     = | \frac{1}{z'-z}\int_z^{z'} \d_y^\alpha d_\eta^\beta
     q(s,y,\eta)\ d s| \\
     \leq C (1+|\eta|)^{-|\beta| + (|\alpha| + |\beta|)/L} 
     \frac{1}{z'-z}\int_z^{z'}(1+ q(z,y,\eta))^{1 - (|\alpha| +
     |\beta|)/L} ds\\
     \leq C (1+|\eta|)^{-|\beta| + (|\alpha| + |\beta|)/L} 
     \left( 1 + \frac{1}{z'-z}\int_z^{z'} 
     q(z,y,\eta) ds \right)^{1 - (|\alpha| + |\beta|)/L}\\
     = C (1+|\eta|)^{-|\beta| + (|\alpha| + |\beta|)/L} 
     (1+ \widehat{q}_{(z',z)}(y,\eta))^{1 - (|\alpha| +
     |\beta|)/L} ,
  \end{multline*}
by Jensen inequality as $t \mapsto - (1+t)^{1 - (|\alpha| +
  |\beta|)/L}$ is convex when $|\alpha| + |\beta| < L$.
\end{proof}
As a consequence of Lemma~\ref{lemma:exp(-Dq)- QL} we have
\begin{lemma}
  Let $q(z,y,\eta) \in \Con^0([0,Z],S^1(\R^{p}\times \R^{r}))$ that
  satisfies Property (\ref{eq:P_L}). Then
  $\widehat{\rho}_\Delta:=\exp[-\Delta \widehat{q}_{(z',z)}(y,\eta)]$ satisfies 
  Property (\ref{eq:Q_L}).
\end{lemma}
The result of Theorem~\ref{theorem:H^s estimate} thus applies to the
modified thin-slab propagator $\widehat{\G}_{(z',z)}$ (Lemma~\ref{lemma:
change of variables} has to be slightly modified).  The proof of
Lemma~\ref{lemma: z' -> Gz'z lipschitz} applies with the aid of
Proposition~\ref{prop: Hs cont any symbol} as
\begin{multline*}
  (\widehat{\G}_{(z^{(2)},z)}-\widehat{\G}_{(z^{(1)},z)})(u_0) (x') =\\
  - \int_{z^{(1)}}^{z^{(2)}}\iint \exp[i \inp{x'-x}{\xi}
    - (z'-z)\widehat{a}_{(z',z)}(x',\xi)]\
  a(z',x',\xi)\ \ u_0(x)\ d x\ \dslash \xi \ d z'.
\end{multline*}
To adapt the proof of Lemma~\ref{lemma: cont G, dz G} we need
\begin{lemma}
  Let $s \in \R$ and $z'',z \in [0,Z]$.  The map $z' \mapsto
  \d_{z'} \widehat{\G}_{(z',z)}$, for $z' \in [z'',z]$, is continuous with
  values in $L(H^{(s+2)}(X),H^{(s)}(X))$, for $z''-z=\Delta$ small
  enough.
\end{lemma}
\begin{proof}
  We choose $\Delta=z''- z$ sufficiently small such that the results
  of Section~\ref{sec:1} apply.  Let $z^{(1)}$, $z^{(2)} \in
  [z,z'']$. Then we have
  \begin{multline*}
   \d_{z'} \widehat{G}_{(z^{(2)},z)}(x',x) - \d_{z'} \widehat{G}_{(z^{(1)},z)}(x',x)\\
   = - \int \exp[i \inp{x'-x}{\xi}] \left( a(z^{(2)},x',\xi ) 
     \exp[-\textstyle{\int}_z^{z^{(2)}} a(s,x',\xi)\ d s]\right.\\
     \left. - a(z^{(1)},x',\xi ) 
     \exp[-\textstyle{\int}_z^{z^{(1)}} a(s,x',\xi)\ d s] \right) d
   \xi\\
   = A_{(z^{(2)},z^{(1)},z)}(x',x) + B_{(z^{(2)},z^{(1)},z)}(x',x),
 \end{multline*}
where 
\begin{multline*}
    A_{(z^{(2)},z^{(1)},z)}(x',x) := 
    - \int \exp[i \inp{x'-x}{\xi}] \ a(z^{(2)},x',\xi ) \\ 
    \left(\exp[-\textstyle{\int}_z^{z^{(2)}}a(s,x',\xi)\ d s]
    - \exp[-\textstyle{\int}_z^{z^{(1)}}a(s,x',\xi)\ d s] \right)\ d \xi,
 \end{multline*}
and 
\begin{multline*}
  B_{(z^{(2)},z^{(1)},z)}(x',x) := - \int \exp[i \inp{x'-x}{\xi}]\\
   (a(z^{(2)},x',\xi ) - a(z^{(1)},x',\xi))\
   \exp[-\textstyle{\int}_z^{z^{(1)}}a(s,x',\xi)\ d s]\ d \xi.
\end{multline*}
We write
\begin{multline*}
    A_{(z^{(2)},z^{(1)},z)}(x',x) =
    \int_{z^{(1)}}^{z^{(2)}} \int \exp[i \inp{x'-x}{\xi}]\\
     a(z^{(2)},x',\xi )\ a(s,x',\xi)\ \exp[-(s-z) \widehat{a}_{(s,z)}(x',\xi)]\
    d s\ d \xi.
\end{multline*}
and for the associated operator, ${\cal A}_{(z^{(2)},z^{(1)},z)}$ we
obtain by Proposition~\ref{prop: Hs cont any symbol} that $\|{\cal
A}_{(z^{(2)},z^{(1)},z)}\|_{(H^{(s+2)},H^{(s)})} \leq C |z^{(2)} -
z^{(1)}|$. For the second term we can apply Proposition~\ref{prop: Hs
cont any symbol} which gives the estimate, for the associated operator,
$\|{\cal B}_{(z^{(2)},z^{(1)},z)}(x',x)\|_{(H^{(s+2)},H^{(s)})} \leq C
\ p(a(z^{(2)},.) - a(z^{(1)},.))$, with $p$ a seminorm in $S^1(X\times
\R^n)$. The continuity of $z \mapsto a(z,.)$ in $S^1(X\times \R^n)$
(Assumption~\ref{assumpt:general assumption})
yields the result. 
\end{proof}
With the previous lemma we can easily adapt the proof of
Lemma~\ref{lemma: cont G, dz G} and obtain the same result for
$\widehat{\G}_{(z',z)}$.
\begin{lemma}
  \label{lemma: cont hat G, dz hat G}
  Let $s \in \R$, $z'',z \in [0,Z]$, with $z <  z''$, and let $u_0
  \in H^{(s+1)}(X)$.  Then the map $z'\mapsto \widehat{\G}_{(z',z)}(u_0)$ is in
  $\Con^0([z,z''],H^{(s+1)}(X)) \cap \Con^1([z,z''],H^{(s)}(X))$ for
  $z''-z=\Delta$ small enough.
\end{lemma}
This allows to use the energy estimate~(\ref{eq:energy estimate}).

We now note that in the proof of Lemma~\ref{lemma: Delta coming out of
symbol}, with the new thin-slab propagator, $\widehat{\G}_{(z',z)}$, the
amplitudes of the operators ${\cal D}_1, \dots,{\cal D}_4$ only
involve the term $a(z',x,\xi)$ instead of both $a(z',x,\xi)$ and
$a(z,x,\xi)$ (as $\d_{z'} ((z'-z)\widehat{a}_{(z',z)}(x',\xi)) =
a(z',x,\xi)$).  Thus the proof of Lemma~\ref{lemma: Delta coming out
of symbol} does not require any assumption like
Assumption~\ref{assumpt:Lipschitz assumption2} made in
Theorem~\ref{theorem: convergence result} or assumptions of
H\"{o}lder type regularity on the symbol $a(z,.)$ made in
Theorem~\ref{theorem: convergence result hoelder}. Consequently we
obtain
\begin{theorem}
  Let $s \in \R$. Then the approximation Ansatz $\widehat{\W}_{\P,z}$
  converges in $L (H^{(s+1)}(\R^n), H^{(s)}(\R^n))$ to the solution
  operator $U(z,0)$ of the Cauchy
  problem~(\ref{eq:one-way})-(\ref{eq:init_cond}) as $\Delta_\P$ goes
  to 0 with a convergence rate of order $\hf$:
  \begin{align*}  
    \| \widehat{\W}_{\P,z}  - U(z,0) \|_{(H^{(s+1)},H^{(s)})} \leq C \Delta_\P^\hf.
  \end{align*}
\end{theorem}

\appendix
\renewcommand{\thefigure}{\arabic{figure}}
\section{ A diagonalization/decoupling of the acoustic wave equation}
\label{app:A}
\newcommand{\matr}[4]{
\left(
\begin{array}{cc}
  #1 & #2\\
  #3 & #4
\end{array}
\right)
}
\newcommand{\vect}[2]{
\left(
\begin{array}{c}
  #1 \\
  #2
\end{array}
\right)
}

We give here an overview of \cite{Stolk:04-bis}, which gives a
motivation for approximating solutions of the Cauchy problem
(\ref{eq:one-way})-(\ref{eq:init_cond}), for instance in the context of
geophysics.

We first consider the scalar wave equation 
\begin{align}
  \label{eq:scalar wave equation}
  \left(  -\rho^{-1} c^{-2} \d_t^2 
    + \sum_{j=1}^{n} \d_j \rho^{-1} \d_j \right) u = F,
\end{align}
as encountered in acoustics, where $\rho$ is the fluid density, and
$c$ is the wavespeed. Both these functions are assumed to be independent
of time $t$ and to be in $\Cinf(\R^n)$.  We further assume that $0 <
\rho_0\leq \rho(y) \leq \rho_1$ and $0< c_0\leq c(y) \leq c_1$, $y \in
\R^n$. We denote $z=y_n$ and $x=(y_1,\dots,y_{n-1})$ and write
$p(x,z,D_t,D_x,D_z) = \rho^{-1} c^{-2} D_t^2 - \sum_{j=1}^{n-1} D_j
\rho^{-1} D_j - D_z \rho^{-1} D_z$ where $D = \frac{1}{i} \d$. Its
principal symbol is $p_2(t,x,z,\tau,\xi,\zeta) =\rho^{-1} (c^{-2}
\tau^2 - |\xi|^2 -\zeta^2)$.
 
Note that $\tau \neq 0$ in $\Char(p)$. We put (\ref{eq:scalar wave
  equation}) in a matrix form
\begin{multline}
  \label{eq: matrix wave equation}
  D_z w(t,x,z) =
 G(x,z,D_t,D_x) w(t,x,z) + f(t,x,z) 
  \mod{\Cinf},\\
  \mbox{with}\ G=\matr{0}{\Lambda \rho}{A}{0}, 
  \quad w=\vect{\Lambda u}{\rho^{-1} D_z u},\quad  f=\vect{0}{F},
\end{multline}
where $\Lambda$ is a first-order elliptic $\psi$DO,
say for instance $|D_{t,x}|$, and
\begin{align*}
  A=\rho^{-1} c^{-2} D_t^2 \Lambda^{-1} 
  - \sum_{j=1}^{n-1} D_j \rho^{-1} D_j \Lambda^{-1},
\end{align*}
with $\Lambda^{-1}$ denoting a parametrix for $\Lambda$.

\begin{figure}

\begin{center}
\includegraphics[width=7cm]{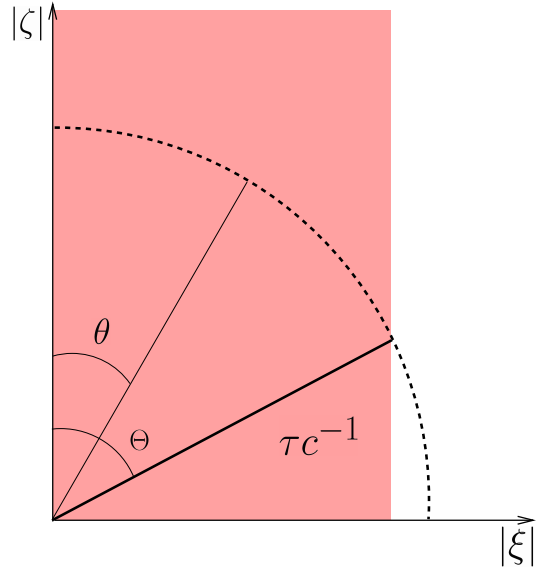}
\caption{The shaded area corresponds to $I_\Theta$ at a given $(t,x,z)$
and a given frequency $\tau$. $\theta$ is the propagation angle. The
set $\Char(p)$ is represented dotted.}
\label{fig:dispersion}
\end{center}
\end{figure}

Following \cite{Stolk:04-bis}, we introduce 
\begin{align*}
  I_\Theta' &= \{ (x,z,\tau,\xi)\ |\ \tau\neq 0 , 
  |c(x,z) \tau^{-1} \xi|\leq \sin \Theta \},\\
  I_\Theta  &= \{ (t,x,z,\tau,\xi,\zeta)\ |\ (x,z,\tau,\xi)\in I_\Theta',\
  |\zeta| \leq c_0^{-1} |\tau|\},
\end{align*}
where $\Theta \in (0,\frac{\pi}{2})$.  The inequality $|\zeta| \leq
c(x,z)^{-1} |\tau|$ on $\Char (p)$ explains the condition $|\zeta|
\leq c_0^{-1} |\tau|$ above. We choose an angle $\Theta\in
(0,\frac{\pi}{2})$ and work in the microlocal region $I_\Theta$
assuming that $\WF(u) \subset I_\Theta$. Figure~\ref{fig:dispersion}
illustrates the set $I_\Theta$ at a given $(x,z)$ and a given frequency
$\tau$.  An angle $\theta \in [-\Theta,\Theta]$ corresponds to a
propagation angle.  Restricting the analysis to $I_\Theta$ corresponds
to staying away from horizontal propagation.  Note that in $I_\Theta$
we have $c(x,z)^{-2} \tau^2 - |\xi|^2>0$, which is the main purpose of
the restriction to such a microlocal region.

In $I_\Theta$, $G$ is a first-order $\psi$DO by
Theorem 18.1.35 in \cite{Hoermander:V3}. In $I_\Theta$ we can follow
the method of \cite[Chapter IX]{Taylor:81} (see also \cite{Taylor:75})
to decouple the up-going and down-going wavefields. We briefly recall
the method here.  Define $\eta^{\pm}(x,z,\tau,\xi) = \pm (c(x,z)^{-2}
\tau^2 - |\xi|^2)^{\hf}$, which are the two roots of $\det (\eta I_2 -
G_1)=0$ with $G_1$ the (matrix-)principal symbol of $G$. The matrix
$G_1(x,z,\tau,\xi)$ is diagonalizable and we choose a matrix
$V(x,z,\tau,\xi) \in S^0(I_\Theta')$,
invertible, such that $V G_1 V^{-1}$ is diagonal; $V$ van be chosen
homogeneous of degree 0. If we write $w^{(0)} = V(x,z,D_t,D_x) w$ we
obtain
\begin{multline*}
  D_z w^{(0)} (t,x,z)= G^{(0)}w^{(0)}(t,x,z) + f^{(0)} (t,x,z)
  \mod{\Cinf}, \\ 
  G^{(0)}= (D_z V) V^{-1}+ V
  G V^{-1} \mod{\Psi^{-\infty}}\ \mbox{in}\ I_\Theta, 
  \quad f^{(0)} = V f.
\end{multline*} 
We write $G^{(0)}=G^{(0)}_1 + G^{(0)}_0$ with $G^{(0)}_1 \in \Psi^1$
in $I_\Theta$ and diagonal and $G^{(0)}_0 \in \Psi^0$ in $I_\Theta$.
By $V^{-1}$ we denote a parametrix for $V$ with principal symbol
$V(x,z,\tau,\xi)^{-1}$ (an abuse of notations, which will occur below again).

We then write $w^{(1)} = (1+K^{(1)}(x,z,D_t,D_x)) w^{(0)}$, with $K^{(1)}
\in \Psi^{-1}$ in $I_\Theta$ of the form
\begin{align*}
  K^{(1)} = \matr{0}{K^{(1)}_{1}}{K^{(1)}_{2}}{0}.
\end{align*}
We then obtain
\begin{multline*}
  D_z w^{(1)} = G^{(0)}_1 w^{(1)}  +  [K^{(1)}, G^{(0)}_1]w^{(1)} 
  + G^{(0)}_0 w^{(1)}  +  f^{(1)} + R^{(1)}w^{(1)} \mod{\Cinf}, \\ 
  R^{(1)} \in \Psi^{-1}\ \mbox{in}\ I_\Theta, \quad
  \quad f^{(1)} = (1+K^{(1)})f^{(0)},
\end{multline*} 
making use of 
\begin{align*}
  (1+K^{(1)}) G^{(0)}_1 (1+K^{(1)})^{-1}=
  G^{(0)}_1 + [K^{(1)}, G^{(0)}_1](1+K^{(1)})^{-1} 
\end{align*}
and the fact that $L (1+K^{(1)})^{-1} - L \in \Psi^{m-1}$ if $L \in
\Psi^{m}$.  Lemma 2.1 in \cite{Taylor:75} shows that $K^{(1)}$ can be
chosen so as to have $[K^{(1)}, G^{(0)}_1]+ G^{(0)}_0$ diagonal up to
an operator in $\Psi^{-1}$ in $I_\Theta$.  The procedure goes on by
choosing $K^{(2)} \in \Psi^{-2}$ in $I_\Theta$ in order to diagonalise
the term of order -1, etc. We thus obtain $Q \in \Psi^{0}$ in
$I_\Theta$ such that $\tilde{w}= Q^{-1}w $ satisfies
\begin{align*}
  D_z \tilde{w}  =  \tilde{G}\tilde{w} +  \tilde{f}  \mod{\Cinf},
  \quad \tilde{f}  = Q^{-1} f,
\end{align*} 
with $\tilde{G} =\tilde{G}(x,z,D_t,D_x) \in \Psi^{0}$ in $I_\Theta$,
diagonal up to a regularizing operator
\begin{align*}
  \tilde{G} = \matr{b_+}{0}{0}{b_-}.
\end{align*}
In \cite{Stolk:04-bis}, Stolk shows that $b_\pm$ can be chosen
selfadjoint.  This is achieved by first choosing selfadjoint operators
with principal symbols equal to $\eta^\pm(x,z,\tau,\xi)$ and then replace
$(1+K^{(i)})$ by $\exp[K^{(i)}]$ in the iteration process described
above. Various choices of $Q$ are presented in \cite{Stolk:04-bis}.

We define the set $J_{\Theta+}$ of points $(t_0,x_0,z_0,\tau_0,\xi_0,
\zeta_0)$ such that the bicharacteristics associated with $b_+$,
parametrized by $z$, $(t(z),x(z),\tau(z),\xi(z))$, passing through
$(t_0,x_0,\tau_0,\xi_0)$ at $z=z_0$, is such that for all $z \in
[0,Z]$, the point $(x(z),z,\xi(z),\tau(z))$ remains in $I_\Theta'$. In
other words, with the interpretation given by
Figure~\ref{fig:dispersion} the propagation angle, $\theta(z)$ along
the bicharacteristics should never exceed $\Theta$.

We now choose $0< \Theta_1 < \Theta_2 < \frac{\pi}{2}$. We choose a
real non-negative symbol $c(z,x,\tau,\xi) \in S^1(\R\times
\R^{n-1}\times \R \times \R^{n-1})$ such that $c=0$ in $I_{\Theta_1}$
and elliptic in the complement of $I_{\Theta_2}$. After extendig
smoothly $b_+$ outside $I_\Theta$, such that  $b_+$ is real
homogeneous of degree 1, we now consider the Cauchy problem
\begin{align*}
(\d_z -i b_+(z,x,D_t,D_x) +c(z,x,D_t,D_x)) v &=0, \\
v(0,.) &= v_+(0,.),
\end{align*}
where 
\begin{align*}
  \tilde{w} = \vect{v_+}{v_-} = Q^{-1} w 
= Q^{-1}\vect{\Lambda u}{\rho^{-1} D_z u}.
\end{align*}
With Assumption (33) and (34) in \cite{Stolk:04-bis}
we obtain that  
\begin{align*}
  &v = v_+ \mod \Cinf\ \mbox{in}\  J_{\Theta_1+},\\
  &v= 0\ \  \mod \Cinf\ \mbox{in the complement of}\ J_{\Theta_2+}.\
\end{align*}
See \cite{Stolk:04-bis} and \cite{Stolk:04} for details.
A similar results holds for the other `one-way' wave operator
$\d_z -i b_- +c$.



\paragraph{Acknowledgement:} 
The author wishes to thank G. H\"ormann for numerous discussions on
many proofs in the paper, especially that of Theorem~\ref{theorem:H^s
estimate} and dicussions on the content of Appendix~\ref{app:A}.

\bibliographystyle{abbrv}

\bibliography{FIO1}

\newcommand{\SortNoop}[1]{}
\begin{thebibliography}{10}

\bibitem{deHoop:04}
M.~V. de~Hoop.
\newblock Microlocal analysis of seismic inverse scattering.
\newblock In G.~Uhlmann, editor, {\em Inside Out, Inverse Problems and
  Applications}, Cambridge, 2004. Cambridge University Press.

\bibitem{dHlRB:03}
M.~V. de~Hoop, J.~H.~L. Rousseau, and B.~Biondi.
\newblock Symplectic structure of wave-equation imaging: A path-integral
  approach based on the double-square-root equation.
\newblock {\em Geoph. J. Int.}, 153:52--74, 2003.

\bibitem{dHlRW:00}
M.~V. de~Hoop, J.~H.~L. Rousseau, and R.-S. Wu.
\newblock Generalization of the phase-screen approximation for the scattering
  of acoustic waves.
\newblock {\em Wave Motion}, 31:43--70, 2000.

\bibitem{Duistermaat:96}
J.~J. Duistermaat.
\newblock {\em Fourier integral operators}.
\newblock Birkh{\"a}user, Boston, 1996.

\bibitem{EN:99}
K.-J. Engel and R.~Nagel.
\newblock {\em One-parameter semigroup for linear evolution equations}.
\newblock Springer-Verlag, Berlin, 1999.

\bibitem{Hoermander:71}
L.~H{\"o}rmander.
\newblock {F}ourier integral operators {I}.
\newblock {\em Acta Math.}, 127:79--183, 1971.

\bibitem{Hoermander:83}
L.~H{\"o}rmander.
\newblock {$L^2$} estimate for {F}ourier integral operators with complex phase.
\newblock {\em Ark. Mat.}, 21(2):283--307, 1983.

\bibitem{Hoermander:V3}
L.~H{\"o}rmander.
\newblock {\em The analysis of linear partial differential operators}, volume
  III.
\newblock Springer-Verlag, 1985.
\newblock Second printing 1994.

\bibitem{Hoermander:V4}
L.~H{\"o}rmander.
\newblock {\em The analysis of linear partial differential operators},
  volume~IV.
\newblock Springer-Verlag, 1985.

\bibitem{Hoermander:V1}
L.~H{\"o}rmander.
\newblock {\em The analysis of linear partial differential operators},
  volume~I.
\newblock Springer-Verlag, second edition, 1990.

\bibitem{K:70}
T.~Kato.
\newblock Linear evolution equations of "hyperbolic" type.
\newblock {\em J. Fac. Sci. Univ. Tokyo, Sec. I}, 17:241--258, 1970.

\bibitem{KiKu:81}
H.~Kitada and H.~Kumano-go.
\newblock A family of {F}ourier integral operators and the fundamental
  solultion for a {S}chr\"odinger equation.
\newblock {\em Osaka J. Math.}, 18:291--360, 1981.

\bibitem{Kumano-go:76}
H.~Kumano-go.
\newblock A calculus of {F}ourier integral operators on {$\R^N$} and the
  fundamental solution for an operator of hyperbolic type.
\newblock {\em Comm. Part. Diff. Eqs.}, 1(1):1--44, 1976.

\bibitem{Kumano-go:81}
H.~Kumano-go.
\newblock {\em Pseudo-differential operators}.
\newblock MIT Press, Cambridge, 1981.

\bibitem{KuTa:79}
H.~Kumano-go and K.~Taniguchi.
\newblock Fourier integral operators of multi-phase and the fundamental
  solution for a hyperbolic system.
\newblock {\em Funkcialaj Ekvacioj}, 22:161--196, 1979.

\bibitem{KuTaTO:78}
H.~Kumano-go, K.~Taniguchi, and Y.~Tozaki.
\newblock Multi-products of phase functions for {F}ourier integral operators
  with an application.
\newblock {\em Comm. Part. Diff. Eqs.}, 3(4):349--380, 1978.

\bibitem{N.Kumano-go:95}
N.~Kumano-go.
\newblock A construction of the fundamental solution for {S}chr\"odinger
  equations.
\newblock {\em J. Math. Sci. Univ. Tokyo}, 2:441--498, 1995.

\bibitem{MeSj:76}
A.~Melin and J.~Sj{\"o}strand.
\newblock Fourier integral operators with complex phase functions and
  parametrix for an interior boundary value problem.
\newblock {\em Comm. Part. Diff. Eqs.}, 1(4):313--400, 1976.

\bibitem{Pazy:83}
A.~Pazy.
\newblock {\em Semigroups of linear operators and applications to partial
  differential equations}.
\newblock Springer-Verlag, New York, 1983.

\bibitem{Ruzhansky:00}
M.~Ruzhansky.
\newblock {\em Regularity theory of {F}ourier integral operators with complex
  phases and singularities of affine fibrations}, volume 131.
\newblock CWI tract, Amsterdam, 2000.

\bibitem{Stolk:04}
C.~C. Stolk.
\newblock Parametrix for a hyperbolic intial value problem with dissipation in
  some region.
\newblock {\em preprint, Ecole Polytechnique, France}, 2004.

\bibitem{Stolk:04-bis}
C.~C. Stolk.
\newblock A pseudodifferential equation with damping for one-way wave
  propagation in inhomogeneous media.
\newblock {\em preprint, Ecole Polytechnique, France}, 2004.

\bibitem{SdH:02}
C.~C. Stolk and M.~V. de~Hoop.
\newblock Microlocal analysis of seismic inverse scattering in anisotropic,
  elastic media.
\newblock {\em Comm.\ Pure Appl.\ Math.}, 55:261--301, 2002.

\bibitem{Taylor:75}
M.~E. Taylor.
\newblock Reflection of singularities of solutions to systems of differential
  equations.
\newblock {\em Commun. Pure Appl. Math.}, 28:457--478, 1975.

\bibitem{Taylor:81}
M.~E. Taylor.
\newblock {\em Pseudodifferential operators}.
\newblock Princeton University Press, Princeton, New Jersey, 1981.

\end{thebibliography}

\end{document}